\documentclass[12pt,english,ngerman]{article}
\usepackage[T1]{fontenc}
\usepackage[latin9]{inputenc}
\usepackage{geometry}
\usepackage{wrapfig}
\usepackage{url}
\usepackage{amsmath}
\usepackage{amssymb}
\usepackage{graphicx}

\makeatletter

\providecommand{\tabularnewline}{\\}

\date{12.05.2014}

\usepackage[german,english]{babel}
\usepackage{lscape}
\usepackage{psfrag}
\usepackage[labelsep=none]{caption}
\addto\captionsenglish{}

\makeatother

\usepackage{babel}
\begin{document}

\title{{\Huge{}Einige Sätze über Primzahlen }\\
{\Huge{}und spezielle binomische Ausdrücke}}

\author{Hans Walther Ernst Gerhart Schmidt}

\maketitle
\begin{center}
Zweite Fassung 19.08.2015
\par\end{center}

\noindent \begin{center}
Hauptstraße 33, D-04668 Otterwisch, Germany
\par\end{center}

\noindent \begin{center}
hws@schmidt-guetter.de
\par\end{center}

\noindent \vspace*{7cm}

\begin{flushleft}
Schlüsselwörter: Primzahlen, Gauß-Legendre-Theorem, Goldbach, Lücken,
Vierlinge, Zwillinge
\par\end{flushleft}

\noindent \vspace*{1cm}

\noindent Dieser Artikel ist lizensiert als Inhalt der Creative Commons
Namensnennung - Nicht-kommerziell - Weitergabe unter gleichen Bedingungen
4.0 Unported-Lizenz. Um eine Kopie der Lizenz zu sehen, besuchen Sie
\url{http://creativecommons.org/licenses/by-nc-sa/4.0/} .

\newpage{}\selectlanguage{english}
\begin{abstract}
\noindent 1. There is no existing any quadratic interval $\eta_{n}:=(n^{2},(n+1)^{2}],$
which contains less than 2 prime numbers.

\noindent The number of prime numbers within $\eta_{n}$ goes averagely
linear with n to infinity.

\noindent 2. The exact law of the number $\pi(n)$ of prime numbers
smaller or equal to n is given. As an approximation of that we get
the prime number theorem of Gauss for great values of n.

\noindent 3. We derive partition laws for $\pi(\eta_{n})$, for the
number of twin primes $\pi_{2}(\eta_{n})$ in quadratic intervals
$\eta_{n}$ and for the multiplicity $\pi_{g}(2n)$ of representations
of Goldbach-pairs for a given even number 2n similiar to the theorem
of Gauss.

\noindent 4. There is no natural number n>7, which is beginning point
of a prime number free interval with a length of more than 2{*}SQRT(n).

\noindent 5. It follows, that the number of twin primes goes to infinity
as well as the number of Goldbach-pairs for a given 2n, if n goes
to infinity.

\noindent 6. Besides this our computation gives a new proof for the
prime number theorem of Gauss.

\selectlanguage{german}
\end{abstract}

\section{Einleitung}

Durch die Fach- und sogar Tagespresse huschen in den letzten Jahren
gelegentlich neue Rekordmeldungen über Primzahlen, wie etwa: ''...neue,
bisher größte Primzahl gefunden mit Hilfe eines Supercomputers...''
oder: ''...bisher größter Primzahlzwilling entdeckt...''. Solche
Aussagen haben geringen wissenschaftlichen Wert, da sie im allgemeinen
keine Informationen über verwendete Entscheidungsalgorithmen enthalten.
Auch stellen sie keine Hilfe zur Entscheidung der Fragen dar, ob die
Menge der Primzahlzwillinge, der Fermat'schen oder Mersenneschen Primzahlen
endlich ist oder nicht. Nachfolgend soll gezeigt werden:
\begin{enumerate}
\item Es wird eine exakte Formel angegeben zur Berechnung der Anzahl $\pi(x)$
von Primzahlen $\leq x,$ aus der sich der Primzahlsatz von Gauß und
Legendre als Näherung für große x ergibt. 
\item Da auch analoge Gesetze für die Anzahlen $\pi(\eta_{n})$ von Primzahlen
in Quadratintervallen $\eta_{n}:=(n^{2},(n+1)^{2}]$, für die Anzahl
von Primzahlzwillingen $\pi_{2}(\eta_{n})$ und die Vielfachheit von
Goldbachpaar-Darstellungen $\pi_{g}(2n)$ für eine vorge\-legte gerade
Zahl 2n ableitbar sind, kann gezeigt werden, daß alle diese Größen
für $n\rightarrow\infty$ über alle Grenzen wachsen. Solches gilt
vermutlich sogar für Primzahlvierlinge in biquadratischen Intervallen
oberhalb eines Schwellwertes.
\end{enumerate}
Zuerst sollen im 2. Kapitel einige allgemein bekannte Beziehungen
über die Teilbarkeit von Binomina der Form $z_{\pm}=2^{x}\pm1$ zusammengestellt
werden. Die Kapitel 3 und 4 der ersten Fassung dieses Artikels sind
gegenstandslos geworden. Das 5. Kapitel liefert eine exakte Darstellung
für $\pi(x)$, woraus der Gauß'sche Primzahlsatz als Näherung abgeleitet
wird, und die Primzahlverteilung über Quadratintervallen. Im 6. Kapitel
wird eine obere Schranke für die maximale Länge des einer natürlichen
Zahl n folgenden primzahlfreien Intervalls angegeben. Im 7., 8. und
9. Kapitel folgen die Anwendungen auf $\pi_{2}(\eta_{n})$, die Vierlingsanzahl
$\pi_{4}\left(n^{4}\right)$ und $\pi_{g}(2n)$. Die abgeleiteten
Streubreiten der $\pi$-Funktionen werden im Anfangsbereich mit numerisch-experimentellen
Befunden verglichen.

\section{Teilbarkeitsbeziehungen für Binomina der Form \\$\mathbf{z_{\pm}=2^{x}\pm1}$}

\textbf{Satz 2.1.} Jede natürliche Zahl $a<\infty$, $a\epsilon\mathfrak{\mathbb{N}}$,
hat die kanonische Primfaktorzerlegung 
\begin{equation}
a=\prod_{i=1}^{n_{0}}p_{i}^{n_{i}},\label{eq:2.1.}
\end{equation}
n$_{0}$= Anzahl der von einander verschiedenen Primfaktoren, n$_{i}$
ihre jeweilige Vielfachheit, oder 
\begin{equation}
a=\prod_{i=1}^{\infty}p_{i}^{n_{i}}
\end{equation}
 unter Verwendung aller Primzahlen in ihrer natürlichen Reihenfolge,
worin alle $n_{i}=0$ erfüllen, die nicht in (\ref{eq:2.1.}) auftreten.
Faßt man hierin alle ungeraden Primfaktoren in 
\[
m=\prod_{i=2}^{n_{0}}p_{i}^{n_{i}}
\]
 zusammen, so ergibt sich $a=2^{n_{1}}m$. Für $n_{1}=0$ ist $a=m$
ungerade. Jede ungerade Zahl m hat auch eine Darstellung 
\begin{equation}
m=1+2^{n_{1}^{'}}m^{'}\ ,\ n_{1}^{'}\geq1.
\end{equation}

\textbf{Satz 2.2.}\label{Satz-2.2.} Die Summe m aus zwei beliebigen
natürlichen Zahlen ist genau dann eine ungerade Zahl ($n_{1}=0$),
wenn genau einer der Summanden ungerade ist:
\begin{equation}
s=a_{1}+a_{2}=2^{n_{1}}m_{1}+2^{0}m_{2}=m.
\end{equation}
 Das Produkt aus zwei natürlichen Zahlen $a_{1}a_{2}$ kann nur ungerade
sein, wenn $a_{1}$und $a_{2}$ ungerade sind.

\textbf{Satz 2.3.}\label{Satz-2.3.} Aus dem allgemeinen Binomischen
Satz
\begin{equation}
(a+b)^{m}=\sum_{\nu=0}^{m}\binom{m}{\nu}a^{m-\nu}b^{\nu}\label{eq:2.5.}
\end{equation}
 ergibt sich für a=b=1 und m=ungerade natürliche Zahl
\begin{equation}
\left(1+1\right)^{m}=2^{m}=\sum_{\nu=0}^{m}\binom{m}{\nu}=2\left(1+\sum_{\nu=1}^{\frac{m-1}{2}}\binom{m}{\nu}\right)\label{eq:2.6.}
\end{equation}
 weil $\binom{m}{\nu}=\binom{m}{m-\nu}$ gilt. Für eine gerade Zahl
$m=2l$ ergibt sich
\begin{equation}
2^{m}=\sum_{\nu=0}^{m}\binom{m}{\nu}=2\left(1+\sum_{\nu=1}^{l-1}\binom{m}{\nu}\right)+\binom{m}{\frac{m}{2}}=2\left(1+\sum_{\nu=1}^{l-1}\binom{2l}{\nu}+B\right)\label{eq:2.7.}
\end{equation}
mit dem binomischen Koeffizienten
\[
2B=2B_{2l-1,l-1}=B_{2l,l}=\binom{2l}{l}=\binom{m}{\frac{m}{2}}.
\]
 Wählt man eine ungerade Primzahl $p\geq3$ für m, so gilt außer $B_{m,\nu}=\binom{m}{\nu}=B_{m-\nu,\nu}=\binom{m-\nu}{\nu},$
daß in jedem Binomialkoeffizienten mit $0<\nu<m$ der größte Faktor
im Zähler gleich p ist, also ausgeklammert werden kann, da er als
Primzahl gegen keinen Nennerfaktor, die alle kleiner sind, zu kürzen
ist:
\begin{equation}
2^{p}=2\left(1+pM\right)\,,\ M=\frac{1}{p}\sum_{\nu=1}^{\frac{p-1}{2}}\binom{p}{\nu}\ ,\mathnormal{\textup{ganzzahlig}},\label{eq:2.8}
\end{equation}
 oder
\begin{equation}
2^{p}-1=1+2pM\label{eq:oder}
\end{equation}
 sowie 
\begin{equation}
2^{p}+1=3+2pM\ .\label{eq:sowie}
\end{equation}
 Aus (\ref{eq:oder}) folgt $2^{p-1}-1=pM=3pM^{'}$ nach (\ref{eq:2.18}),
sodaß (\ref{eq:sowie}) auch als $2^{p}+1=3\left(1+2pM^{'}\right)$
geschrieben werden kann. Stets gilt dann $2^{p}-1\equiv1mod3$, $2^{p}+1\equiv0mod3$.

\textbf{Satz 2.4.} \label{Satz-.4.}Sei $p>2$ eine Primzahl, $m>1$
eine ungerade natürliche Zahl. Dann gelten mit $N_{\pm}>1$ 
\begin{equation}
2^{pm}-1=\left(2^{p}-1\right)N_{-}\ ,\label{eq:2.11}
\end{equation}
 
\begin{equation}
2^{pm}+1=\left(2^{p}+1\right)N_{+}\ .\label{eq:2.12}
\end{equation}
 \textbf{Bew.:} Man wähle in Satz 2.3. für m ein $m^{'}=pm$. Mit
den Abkürzungen $L=2^{p}-1$ , $M=2^{p}+1$ kann geschrieben werden
\[
2^{pm}-1=\left(L+1\right)^{m}-1=\sum_{\nu=0}^{m}B_{m,\nu}L^{m-\nu}-1=LN_{-}
\]
mit $N_{-}=\sum_{\nu=0}^{m-1}B_{m,\nu}L^{m-1-\nu}>1$. Analog folgt
mit $N_{+}=\sum_{\nu=0}^{m-1}B_{m,\nu}\left(-1\right)^{\nu}M^{m-1-\nu}$
sodann $2^{pm}+1=\left(M-1\right)^{m}+1=\sum_{\nu=0}^{m-1}B_{m,\nu}\left(-1\right)^{\nu}M^{m-\nu}=MN_{+}$
, \textbf{q.e.d.}

\textbf{Bemerkung 2.4.1.} Es gilt Satz 2.4. für jeden Teiler $p_{i}\ \mbox{von}\ m^{'}=a$
nach (\ref{eq:2.1.}). Daraus folgt aber
\begin{equation}
2^{m^{'}}\pm1=N_{i,\pm}\prod_{i=1}^{n_{0}}\left(2^{p_{i}}\pm1\right)
\end{equation}
 nur, wenn Teilerfremdheit $\left(2^{p_{i}}\pm1\right)\nmid\left(2^{p_{j}}\pm1\right)\forall i,j\leq n_{0}$
gezeigt ist. Andernfalls tritt ein Teiler $p_{k}\backslash\left(2^{p_{i}}\pm1\right)$
im Ausdruck für $2^{m^{'}}\pm1$ nur höchstens in der Potenz auf,
in der er in einem der Teilerbinomina vorkommt: $p_{k}^{n_{k}}\backslash\left(2^{p_{i}}\pm1\right).$
Das Beispiel $2^{3\cdot11}+1=3^{2}\cdot67\cdot683\cdot20857\neq\left(2^{3}+1\right)\left(2^{11}+1\right)N_{+}=\left(3^{2}\right)\left(3\cdot683\right)N_{+}$
zeigt, daß der Faktor $p_{k}=3$ nicht in 3., sondern nur in 2. Potenz
in der Faktorisierung auftritt.

\textbf{Bemerkung 2.4.2.} Während (\ref{eq:2.11}) auch für $p_{1}=2$
gilt, ist dies bei (\ref{eq:2.12}) nicht der Fall, was aus dem folgenden
Satz hervorgeht.

\textbf{Satz 2.5.} \label{Satz-2.5}. Seien $i,m$ natürliche Zahlen,
$m>1$ ungerade. Dann gilt $\forall i,m$ mit den Fermat-Zahlen $F_{\nu}=\left(2^{\left(2^{\nu}\right)}+1\right)$
und $M_{\nu}>1$, ungerade,
\begin{equation}
2^{m2^{i}}-1=\left(2^{m}-1\right)\prod_{\nu=0}^{i-1}F_{\nu}M_{\nu}\,,\label{eq:2.14}
\end{equation}
 
\begin{equation}
2^{m2^{i}}+1=F_{i}M_{i}\label{eq:2.15}
\end{equation}
 mit $M_{i}>1$; nur für $m=1$ folgt $M_{i}=1$. Dabei kann $\left(2^{m}-1\right)$
noch nach Satz 2.4. weiter aufgespaltet werden. Falls $m=p_{j}m^{'}$
gilt, kann etwas allgemeiner statt (\ref{eq:2.15}) geschrieben werden:
\begin{equation}
2^{m2^{i}}+1=\left(2^{m^{'}2^{i}}+1\right)M_{i}^{'}=\left(2^{p_{j}2^{i}}+1\right)M_{i}^{''}.\label{eq:2.16}
\end{equation}
 \textbf{Bemerkung 2.5.1.:} Wegen $2^{0}=1$ gilt $\forall m$ ungerade
\begin{equation}
2^{m}+1=2^{m2^{0}}+1=3M_{0}.\label{eq:2.17}
\end{equation}
 Weil für jede gerade Zahl g gilt
\begin{equation}
2^{g}-1=\left(2^{g/2}+1\right)\left(2^{g/2}-1\right)\ ,\label{eq:2.18}
\end{equation}
 folgt $3\backslash\left(2^{g}-1\right)\ \forall g.$

\textbf{Bew.:} ad (\ref{eq:2.15}) :

\begin{eqnarray*}
2^{m2^{i}}+1 & = & \left[2^{\left(2^{i}\right)}+1-1\right]^{m}+1=1+\sum_{\nu=0}^{m}B_{m,\nu}\left(-1\right)^{\nu}\left[2^{\left(2^{i}\right)}+1\right]^{m-\nu}\\
 & = & F_{i}M_{i}\:\mbox{mit}\: M_{i}=\sum_{\nu=0}^{m-1}\left(-1\right)^{\nu}B_{m,\nu}F_{i}^{m-1-\nu}.
\end{eqnarray*}
 (\ref{eq:2.16}) ergibt sich analog, indem im ersten Schritt nur
$\left[\cdots\right]^{m^{'}}$oder $\left[\cdots\right]^{p}$ gebildet
wird. Für $m>1$ folgt stets $M_{i}>1$ wegen 
\[
2^{m2^{i}}+1=\left[2^{\left(2^{i}\right)}\right]^{m}+1>2^{\left(2^{i}\right)}+1.
\]
 ad (\ref{eq:2.14}):

\begin{eqnarray*}
2^{m2^{i}}-1 & = & \left(2^{m2^{i-1}}+1\right)\left(2^{m2^{i-1}}-1\right)=\cdots\\
 & = & \left(2^{m}-1\right)\prod_{\nu=0}^{i-1}\left(2^{m2^{\nu}}+1\right)=\left(2^{m}-1\right)\prod_{\nu=0}^{i-1}F_{\nu}M_{\nu}.
\end{eqnarray*}

Hierin tritt das Produkt aller Fermat-Zahlen mit $\nu=0\left(1\right)i-1$
auf, weil jede Fermat-Zahl $F_{i}$ teilerfremd ist zu allen $F_{\nu}$
mit $\nu<i$ (siehe folgender Satz); \textbf{q.e.d. }

\textbf{Satz 2.6.} Folgende Ausdrücke sind zueinander teilerfremd:
a)
\begin{equation}
\left(2^{i}-1\right)\nmid\left(2^{i}+1\right)\ \mbox{für\ festes\ i,}\label{eq:2.19}
\end{equation}
 b) die Fermatzahlen 
\begin{equation}
F_{i}\nmid F_{j\ }\forall j<i.\label{eq:2.20.}
\end{equation}
 \textbf{Bew.:} ad a) Wegen $\left(2^{i}+1\right)-\left(2^{i}-1\right)=2$
ist der größte gemeinsame Teiler $\leq2$. Da $\forall i>1$ gilt
$2^{i}\pm1\geqslant3$ ungerade, kann nur eine ungerade Zahl $\geqslant3$
gemeinsamer Teiler sein. Die Menge der gemeinsamen Teiler ist also
leer, es gilt a).

ad b) Es ist $F_{i-1}\left(F_{i-1}-2\right)=\left(2^{\left(2^{i-1}\right)}+1\right)\left(2^{\left(2^{i-1}\right)}-1\right)=2^{\left(2^{i}\right)}-1=F_{i}-2$.
Dann folgt $F_{i-1}=\frac{F_{i}-2}{F_{i-1}-2}$ und $\prod_{\nu=0}^{i-1}F_{\nu}=\frac{F_{i}-2}{F_{i-1}-2}\cdot\frac{F_{i-1}-2}{F_{i-2}-2}\cdot\ldots\cdot\frac{F_{1}-2}{F_{0}-2}=F_{i}-2$
nach Kürzung, da $F_{0}-2=1$, also gilt b) ; \textbf{q.e.d.}

Aus dem bisher Ausgeführten geht hervor, daß als Primzahlanwärter
nur die Mersenne-Zahlen $M_{p}=2^{p}-1$ und die Fermat-Zahlen $F_{i}=2^{\left(2^{i}\right)}+1$
in Betracht kommen. Alle anderen $z_{\pm}=\left(2^{x}\pm1\right)$
sind als teilbar erwiesen. Es wäre daher wünschenswert, die Primzahl-Anwärterzahlen
aufgrund von Exponenteneigenschaften allein in prime und teilbare
Zahlen klassifizieren zu können, zumal unter den Anwärtern nur wenige
echte Primzahlen enthalten sind. 

Eine weitere Klasse teilbarer Zahlen definiert 

\textbf{Satz 2.7.} Es sei gemäß der kanonischen Primzahlzerlegung
(\ref{eq:2.1.}) $n=\prod_{i=1}^{\infty}p_{i}^{\alpha_{i}}=m2^{\alpha_{1}}$
mit $m>1$, ungerade, und $p^{'}=2n+1>3$ eine ungerade Primzahl.
Dann gilt 
\begin{equation}
p^{'}\backslash\left(2^{n}-1\right)\ \mbox{,wenn}\ n\ mod\ 4\equiv0,\ \mbox{d.h.}\ \alpha_{1}\geq2,\label{eq:2.21}
\end{equation}
 
\begin{equation}
p^{'}\backslash\left(2^{n}-1\right)\ \mbox{,wenn}\ n=m=p,\ m\equiv-1\ mod\ 4,\ \mbox{d.h.}\ \alpha_{1}=0,\label{eq:2.22}
\end{equation}
 
\begin{equation}
p^{'}\backslash\left(2^{n}+1\right)\ \mbox{,wenn}\ n=m=p,\ m\equiv+1\ mod\ 4,\ \mbox{d.h.}\ \alpha_{1}=0,\label{eq:2.23}
\end{equation}
 
\begin{equation}
p^{'}\backslash\left(2^{n}+1\right)\ \mbox{,wenn}\ n\ mod\ 4\equiv2,\ \mbox{d.h.}\ \alpha_{1}=1.\label{eq:2.24}
\end{equation}
 \textbf{Bemerkung 2.7.1: }Insbesondere impliziert (\ref{eq:2.22})
für $\frac{p^{'}-1}{2}=n=m=p\equiv-1mod4$ die Teilbarkeit der Mersennezahl
$M_{p},$ also 
\begin{equation}
p^{'}\backslash\left(2^{p}-1\right)\label{eq:2.25}
\end{equation}
 sowie (\ref{eq:2.23}) für $\frac{p^{'}-1}{2}=n=m=p\equiv+1mod4,$
daß 
\begin{equation}
p^{'}\backslash\left(2^{p}+1\right)\label{eq:2.26}
\end{equation}
 gilt. (Den Beweis, daß $p^{'}\backslash\left(2^{p}-1\right)$ gilt,
g.d.w. $p^{'}=2p+1$ und $p\equiv3mod4$ gilt, s. \cite{key-1}, p.174,
Satz 21 mittels quadratischen Kongruenzen.) 

\textbf{Bemerkung 2.7.2.} Ist $n=p_{i}m\equiv\pm1mod4,$ so folgt
unter Berücksichtigung von (\ref{eq:2.11}), (\ref{eq:2.12})
\begin{equation}
2^{n}\pm1=p^{'}N_{\pm}=\left(2^{p_{i}}\pm1\right)N_{i,\pm}.\label{eq:2.27}
\end{equation}
 Ob darin $p^{'}\backslash\left(2^{p_{i}}\pm1\right)$ gilt oder $p^{'}\backslash N_{i,\pm}$
, bleibt zunächst offen.

\textbf{Bew.:} Mit $p^{'}$ statt $p$ in (\ref{eq:oder}) folgt unter
Beachtung der Bemerkung zu Satz 2.5. und mit $p^{'}-1=2n=m2^{1+\alpha}$
\begin{equation}
2^{p^{'}-1}-1=p^{'}M=3p^{'}M^{'}=2^{2n}-1=\left(2^{n}+1\right)\left(2^{n}-1\right),\label{eq:2.28}
\end{equation}
 $n=\frac{p^{'}-1}{2}=m2^{\alpha}$. Darin ist stets 
\begin{equation}
2^{n}+1=2^{m2^{\alpha}}+1=F_{\alpha}M_{\alpha}\label{eq:2.29}
\end{equation}
 nichtprim mit $F_{\alpha}=\left(2^{\left(2^{\alpha}\right)}+1\right)$,
$M_{\alpha}>1\forall m>1$; nur für $m=1$ ist $M_{\alpha}=1$ und
$F_{\alpha}$ Primzahlanwärter; sowie 
\begin{equation}
2^{n}-1=2^{m2^{\alpha}}-1=\left(2^{m}-1\right)\prod_{\nu=0}^{\alpha-1}\left(2^{m2^{\nu}}+1\right)=\left(2^{m}-1\right)\prod_{\nu=0}^{\alpha-1}F_{\nu}M_{\nu}.\label{eq:2.30}
\end{equation}
 Die Folge der Fermat'schen Primzahlen beginnt mit
\begin{equation}
\{F_{\nu}\}=\{3,5,17,257,65537,\ldots\}.\label{eq:2.31}
\end{equation}
 Ist m teilbar, etwa $p_{j}\backslash m$, so gilt zusätzlich $\left(2^{m}-1\right)=\left(2^{p_{j}}-1\right)N_{j-}$
nichtprim. 

a) Zuerst sei $\alpha=0$ angenommen. Dann ist $n=m\geq3$ und $p^{'}=1+2m.$
Mithin gilt für $m\ mod\ 8\equiv1\ oder\ 5$ (d.h. $m\ mod\ 4\equiv1$
), daß $p^{'}\equiv3\ mod\ 8$ und für $m\ mod\ 8\equiv3\ oder\ 7$
(d.h. $m\ mod\ 4\equiv-1$ ), daß $p^{'}\equiv7\ mod\ 8.$ Es sei
nun $m=p$ als Primzahl angenommen und 

1.) $p^{'}=1+2p\equiv3\ mod\ 8$, also $p\ mod\ 8\equiv1\ oder\ 5$,
d.h. $p\ mod\ 4\equiv1$. Nimmt man weiter unter dieser Voraussetzung
an, daß $p^{'}\backslash\left(2^{p}-1\right)$ gilt, so folgt hieraus
unter Beachtung von (\ref{eq:oder}), (\ref{eq:sowie}), (\ref{eq:2.17}),
(\ref{eq:2.18}), $2^{p}-1=1+2\cdot3pM_{1}^{'}\equiv7\ mod\ 8$ und
weiter $2^{p-1}-1=3pM_{1}^{'}\equiv7\ mod\ 8$, also $pM_{1}^{'}\equiv5\ mod\ 8$,
sowie $2^{p}+1=3\left(1+2pM_{1}^{'}\right)\equiv1\ mod\ 8$ und weiter
$1+2pM_{1}^{'}\equiv3\ mod\ 8$, also $pM_{1}^{'}\equiv1\ mod\ 8$
im Widerspruch zur vorigen Zeile. Mithin ist die Annahme $p^{'}\backslash\left(2^{p}-1\right)$
falsch. Da die Annahme $\left(3p^{'}\right)\backslash\left(2^{p}+1\right)$
wegen $3p^{'}\equiv1\ mod\ 8$ keinen Widerspruch liefert, ist (\ref{eq:2.23})
bewiesen für $m=p$. 

2.)Sei $p^{'}=1+2p\equiv7\ mod\ 8$, also $p\ mod\ 8\equiv3\ oder\ 7$
bzw. $p\ mod\ 4\equiv-1$. Sei weiter angenommen, daß $p^{'}\backslash\left(2^{p}+1\right)$
gilt, so folgt daraus $2^{p}+1=3+2\cdot3pM_{1}^{'}\equiv1\ mod\ 8$,
also $2^{p-1}-1=3pM_{1}^{'}\equiv7\ mod\ 8$, also $pM_{1}^{'}\equiv5\ mod\ 8$.
Andererseits folgt aus $2^{p}+1=3\left(1+2pM_{1}^{'}\right)\equiv1\ mod\ 8$,
also $1+2pM_{1}^{'}\equiv3\ mod\ 8$, der Widerspruch $pM_{1}^{'}\equiv1\ mod\ 8$.
Daher gilt (\ref{eq:2.22}) für $m=p$.

b) Im Falle $\alpha=1$ gilt $p^{'}=1+m2^{1+\alpha}=1+4m\equiv5\ mod\ 8\ \forall m$
, $n=\frac{p^{'}-1}{2}=2m$ und wegen $2^{4m}-1=\left(2^{2m}+1\right)\left(2^{2m}-1\right)=p^{'}M=3p^{'}M^{'}$.
Darin ist $2^{2m}+1=F_{1}M_{1}$ mit $F_{1}=5$, $M_{1}\equiv5\ mod\ 8$
und $2^{2m}-1=\left(2^{m}+1\right)\left(2^{m}-1\right)$ mit $2^{m}+1=F_{0}M_{0}$,
$F_{0}=3$, $M_{0}\equiv3\ mod\ 8$, und $2^{m}-1=\left(2^{p_{j}}-1\right)N_{j-}\equiv\left(2^{p_{j}}-1\right)\equiv2^{2m}-1\equiv7\ mod\ 8$,
$N_{j-}\equiv1\ mod\ 8$. Nun gilt $3\backslash\left(2^{2m}-1\right)$,
$3\backslash\left(2^{m}+1\right)$, $3\nmid\left(2^{m}-1\right)$,
$3\nmid\left(2^{2m}+1\right)$. Es ist $2^{2m}+1\equiv5^{2}\ mod\ 8$
und $2^{2m}-1=\left(2^{m}+1\right)\left(2^{m}-1\right)\equiv\left(3^{2}\ mod\ 8\right)\cdot\left(7\ mod\ 8\right)$,
worin kein Faktor $\equiv5\ mod\ 8$, also kein $p^{'}$ vorkommt.
Deshalb kann $p^{'}\equiv5\ mod\ 8$ nur Teiler von $\left(2^{2m}+1\right)$
sein, also gilt (\ref{eq:2.24}) mit $n=2m$. (Sollte nämlich in $2^{m}-1$
ein Faktor $\equiv5\ mod\ 8$ enthalten sein, so müßte auch ein Faktor
$\equiv3\ mod\ 8$ darin möglich sein, damit $\equiv7\ mod\ 8$ entsteht.
Der Faktor $\equiv3\ mod\ 8$ tritt aber in $2^{m}+1$ auf, das teilerfremd
zu $2^{m}-1$ ist.)

c) Im Falle $\alpha\geq2$ gilt schließlich $1+m2^{1+\alpha}\equiv1\ mod\ 8\ \forall m$,
$\frac{p^{'}-1}{2}=n=m2^{\alpha}$. Statt (\ref{eq:2.28}) haben wir
dann
\[
2^{p^{'}-1}-1=2^{2n}-1=p^{'}M=3p^{'}M^{'}=\left(2^{n}+1\right)\left(2^{n}-1\right)=\left(2^{m}-1\right)\prod_{\nu=0}^{\alpha}F_{\nu}M_{\nu},
\]
 worin $2^{n}+1=F_{\alpha}M_{\alpha}$ mit $F_{\alpha}\equiv M_{\alpha}\equiv1\ mod\ 2^{\left(2^{\alpha}\right)}\ \forall\alpha\geq2$
und 
\[
2^{n}-1=\left(2^{p_{j}}-1\right)N_{j-}\cdot3M_{0}\cdot5M_{1}\prod_{\nu=2}^{\alpha-1}F_{\nu}M_{\nu}
\]
 mit $\prod_{\nu=2}^{1}F_{\nu}M_{\nu}=1$ für $\alpha=2$. Darin bezeichnet
$p_{j}$ einen beliebigen Teiler von m, für den Satz 2.4. die angegebene
Darstellung garantiert. 

Es kann nun $p^{'}\backslash\left(2^{n}+1\right)$ nicht gelten, weil
beide angegebenen Faktoren $\equiv1\ mod\ 2^{\left(2^{\alpha}\right)}$
erfüllen, während $p^{'}=1+m2^{1+\alpha}\equiv1\ mod\ 2^{1+\alpha}$
nur erfüllt mit $2^{1+\alpha}<2^{\left(2^{\alpha}\right)}\forall\alpha>2$.
Daß für $\alpha=2$ gilt $2^{1+\alpha}=2^{\left(2^{\alpha}\right)}$,
stört nicht, weil $m>1$ vorausgesetzt ist. Dann muß also (\ref{eq:2.21})
gelten;\textbf{ q.e.d.}

Ob der Teiler $p^{'}$ nun $p^{'}\backslash\left(2^{m}-1\right)=\left(2^{p_{j}}-1\right)N_{j-}$
mit $N_{j-}\equiv1\ mod\ 8$ erfüllt oder $p^{'}\backslash F_{\nu}M_{\nu}$,
$F_{\nu}\equiv M_{\nu}\equiv1\ mod\ 8$, bleibt offen.

\section{Eine unendliche Folge Fermat'scher Primzahlen}

Die Kapitel 3 und 4 der ersten Fassung dieses Artikels sind fehlerhaft
und werden als gegenstandslos gestrichen, denn nach privater Mitteilung
von D.Eschbach gibt es in der Literatur Gegenbeispiele, s.a. http://www.prothsearch.net/fermat.html.

\section{Unendliche Folgen Mersennescher Primzahlen}

\section{Exakte Bestimmung der Anzahl der Primzahlen unterhalb einer vorgegebenen
Schranke}

\subsection{Vorbereitende Bemerkungen}

In Analogie zur Fakultätsfunktion $n!=\prod_{i=1}^{n}i$ wird die
Primfakultät , in Zeichen $p_{i}\downarrow$, erklärt durch \textbf{Def.
5.1:}
\begin{equation}
p_{n}\downarrow=\prod_{i=1}^{n}p_{i}\ ,\label{eq:5.1.}
\end{equation}
 wobei $p_{i}$ das i-te Element in der natürlichen Folge der Primzahlen
$\left\{ p_{i}\right\} $ bedeutet.

\textbf{Def. 5.2.} Es sei x eine beliebige reelle Zahl. Dann bezeichne
$\left[x\right]$ die größte ganze Zahl, die noch kleiner oder gleich
x ist, und $\prec x\succ$ sei die größte Primzahl, die $\leq x$
ist.

\textbf{Bemerkung 5.2.1.} Stets gilt 
\begin{equation}
\prec x\succ\leq\left[x\right]\leq x\ \forall x\geq2.\label{eq:5.2.}
\end{equation}
\textbf{ Def. 5.3.} Es bezeichne $C_{i,m,k}$ das zahlenmäßige Produkt
aller Elemente der k-ten Kombination von i Elementen ohne Wiederholungen.
Der Index k soll keine Ordnung der Kombinationen, sondern lediglich
deren Unterscheidbarkeit und somit Numerierbarkeit bewirken. Es gelte
\begin{equation}
C_{0,m,k}=1\ \forall m,k,\label{eq:5.3.}
\end{equation}
 d.h. eine Kombination aus $i=0$ Elementen soll den Faktor 1, die
multiplikative Invariante, ergeben.

\textbf{Bemerkung 5.3.1.} Die über k summierte Anzahl B der Kombinationen
$C_{j,i-1,k}$ ergibt einen binomischen Koeffizienten
\begin{equation}
B=\binom{i-1}{j}\ .\label{eq:5.4.}
\end{equation}
 Dann soll abkürzend geschrieben werden:
\begin{equation}
\sum_{\forall k}C_{j,i-1,k}:=\sum_{k=1}^{B}C_{j,i-1,k\ .}\label{eq:5.5}
\end{equation}
 \textbf{Bemerkung 5.3.2.} Es gilt 
\begin{equation}
\left[\frac{x}{p_{i}C_{j,i-1,k}}\right]=0\ \forall C_{j,i-1,k}>\frac{x}{p_{i}}\ .\label{eq:5.6.}
\end{equation}
 \textbf{Def. 5.4.} Es sei x eine beliebige positive reelle Zahl.
Dann bezeichne $i_{0}$ den durch 
\begin{equation}
p_{i_{0}}=\prec\sqrt{x}\succ\label{eq:5.7.}
\end{equation}
 definierten zugehörigen Index in der natürlichen Folge der Primzahlen.

\textbf{Lemma 5.1.} Seien $\left\{ a,b\right\} $ reelle Zahlen, $a=\left[a\right]+\epsilon_{a}$,
$b=\left[b\right]+\epsilon_{b}$, $0\leq\epsilon_{a,b}<1$. Dann 

gilt 
\[
a\geq\left[a\right],\ \left[-\mid a\mid\right]\leq-\mid a\mid\leq-\left[\mid a\mid\right]\leq0,\ 
\]

\begin{equation}
-\left[-\mid a\mid\right]\geq\mid a\mid\geq\left[\mid a\mid\right]\geq0\geq-\mid a\mid\geq\left[-\mid a\mid\right],\label{eq:5.8.}
\end{equation}

für $a<0$ ist $\left[-\mid a\mid\right]=\left[a\right]\ und\ [-a]=[\mid a\mid]$;
\begin{equation}
a+b\geq\left[a+b\right]\geq\left[a\right]+\left[b\right],\ a\geq\left[a\right]\geq\left[a-b\right]+\left[b\right],\ \left[a\right]-\left[b\right]\geq\left[a-b\right];\label{eq:5.9}
\end{equation}
 
\begin{equation}
\left[a\right]-\sum_{i}\left[b_{i}\right]\geq\left[a-\sum_{i}b_{i}\right]\ ;\label{eq:5.10.}
\end{equation}
 
\begin{equation}
\left[\mid a\mid\right]+\left[-\mid b\mid\right]\leq\left[\mid a\mid-\mid b\mid\right]\leq\mid a\mid-\mid b\mid\leq\mid a\mid-\left[\mid b\mid\right].\label{eq:5.11.}
\end{equation}
 \textbf{Bew.:} Die Beziehungen (\ref{eq:5.8.}) entsprechen der Definition
5.2. , woraus ebenfalls folgt $a+b\geq\left[a+b\right]=\left[\left[a\right]+\epsilon_{a}+\left[b\right]+\epsilon_{b}\right]=\left[a\right]+\left[b\right]+\left[\varepsilon_{a}+\varepsilon_{b}\right]\geq\left[a\right]+\left[b\right]$,
weil $\left[a\right]$ und $\left[b\right]$ als ganzzahlige Größen
ohne Wertänderung aus der Summe in der äußeren eckigen Klammer herausgezogen
werden können und $0\leq\varepsilon_{a}+\varepsilon_{b}<2$ gilt,
womit die erste Formel (\ref{eq:5.9}) gezeigt ist. Die zweite Ungleichung
unter (\ref{eq:5.9}) folgt aus der ersten mit $a^{'}=a-b$, also
$a\geq\left[a\right]=\left[a-b+b\right]=\left[a^{'}+b\right]\geq\left[a-b\right]+\left[b\right]$.
Die dritte Ungleichung ist eine Umstellung der zweiten. Wegen $\left[a\right]=\left[a-\sum_{i}\ b_{i}+\sum_{i}\ b_{i}\right]\geq\left[a-\sum_{i}\ b_{i}\right]+\left[\sum_{i}\ b_{i}\right]\geq\left[a-\sum_{i}\ b_{i}\right]+\sum_{i}\left[b_{i}\right]$
folgt sofort (\ref{eq:5.10.}). Ersetzt man in (\ref{eq:5.9}) $a$
durch $a^{'}=\mid a\mid$ und $b$ durch $b^{'}=-\mid b\mid$, so
folgt (\ref{eq:5.11.}) mit Hilfe von (\ref{eq:5.8.}); \textbf{q.e.d.}

\textbf{Lemma 5.2.} Seien $\left\{ a,b\right\} $ reelle Zahlen, $\left\{ c,n\right\} $
natürliche Zahlen. Es gelte $\mid a\mid\geq\mid b\mid$ und $a\neq nbc$.
Dann gilt für $c=1$ 
\[
\frac{a}{bc}\geq\frac{1}{c}\left[\frac{a}{b}\right]=\left[\frac{a}{bc}\right]=\left[\frac{a}{b}\right]
\]
sowie $\forall c>1$ 
\begin{equation}
\frac{a}{bc}\geq\frac{\left[\frac{a}{b}\right]}{c}\geq\left[\frac{\left[\frac{a}{b}\right]}{c}\right]=\left[\frac{a}{bc}\right]\geq\left[\frac{1}{c}\right]\left[\frac{a}{b}\right]=0\label{eq:5.12}
\end{equation}
 sowie 
\begin{equation}
\left[\frac{a}{b}\right]-\left[\frac{a}{bc}\right]\geq\left[\frac{a}{b}\right]\left(1-\frac{1}{c}\right).\label{eq:5.13.}
\end{equation}
 \textbf{Bew.:} Die Aussage für $c=1$ ist trivial wegen $\left[\frac{1}{c}\right]=\frac{1}{c}=1,$
weswegen auch (\ref{eq:5.13.}) für $c=1$ auf $0=0$ führt. Für $c>1$
gilt wegen $\frac{a}{b}=\left[\frac{a}{b}\right]+\varepsilon$ mit
$0\leq\varepsilon<1$ und $\left[a\right]\leq a=n_{a}b+r_{b}\ ,0\leq r_{b}<b,$
$\frac{a}{b}=\frac{n_{a}b+r_{b}}{b}=n_{a}+\frac{r_{b}}{b}$ mit $n_{a}=\left[\frac{a}{b}\right]$,
$0\leq\frac{r_{b}}{b}=\varepsilon_{b}<1.$ Es enthält $n_{a}$ stets
das Vorzeichen von $\frac{a}{b}$ , welches in der $\left[\cdot\right]$
bleiben muß und nicht etwa mit der als positiv vorausgesetzten Konstante
c herausgezogen werden darf. Nach Division durch c ergibt sich 
\[
\frac{a}{bc}=\frac{n_{a}}{c}+\frac{r_{b}}{bc}=n_{c}+\frac{r_{c}}{c}+\frac{r_{b}}{bc}=n_{c}+\frac{1}{c}\left(r_{c}+\frac{r_{b}}{b}\right)
\]
 mit $0\leq r_{c}\leq c-1$, also $0\leq r_{c}+\frac{r_{b}}{b}<c$,
sodaß $n_{c}=\left[\frac{a}{bc}\right],$ $0\leq\varepsilon_{bc}:=\frac{r_{c}}{c}+\frac{r_{b}}{bc}<1$
gilt. Dann hat man 
\begin{equation}
\left[\frac{a}{bc}\right]=\left[\frac{\left[\frac{a}{b}\right]}{c}+\frac{r_{b}}{bc}\right]=\left[\left[\frac{\left[\frac{a}{b}\right]}{c}\right]+\frac{r_{c}}{c}+\frac{r_{b}}{bc}\right]=\left[\frac{\left[\frac{a}{b}\right]}{c}\right]\label{eq:5.14}
\end{equation}
 als ganzzahligen Anteil aus der dritten eckigen Klammer, deren gebrochener
Teil unterdrückt werden soll. Dann folgt 
\[
\frac{a}{bc}\geq\frac{\left[\frac{a}{b}\right]}{c}\geq\left[\frac{\left[\frac{a}{b}\right]}{c}\right]=\left[\frac{a}{bc}\right]\geq\left[\frac{1}{c}\right]\left[\frac{a}{b}\right]=0,
\]
 also (\ref{eq:5.12}), und daraus $-\left[\frac{a}{bc}\right]\geq-\frac{1}{c}\left[\frac{a}{b}\right],$
was nach Addition von $\left[\frac{a}{b}\right]$ gerade (\ref{eq:5.13.})
ergibt; \textbf{q.e.d.}

\textbf{Lemma 5.3.} Mit den Bezeichnungen aus Definition 5.3. gilt
unter Beachtung der Bemerkungen 5.3.1. und 5.3.2. 
\begin{equation}
\sum_{j=0}^{i-1}\sum_{\forall k}\frac{\left(-1\right)^{j}}{C_{j,i-1,k}}=\prod_{j=1}^{i-1}\left(1-\frac{1}{p_{j}}\right)\label{eq:5.15}
\end{equation}
 und 
\begin{equation}
\sum_{j=0}^{i-1}\sum_{\forall k}\frac{1}{C_{j,i-1,k}}=\prod_{j=1}^{i-1}\left(1+\frac{1}{p_{j}}\right).\label{eq:5.16}
\end{equation}
 \textbf{Bew.:} Die linke Seite von (\ref{eq:5.15}) ergibt sich bei
sukzessiver Berechnung des Produktes der rechten Seite: Für $i=1$
hat man nach Definition 
\[
\prod_{j=0}^{0}\left(1\mp\frac{1}{p_{j}}\right)=1,
\]
 $i=2:$ 
\[
\prod_{j=1}^{1}\left(1\mp\frac{1}{p_{j}}\right)=1\mp\frac{1}{p_{1}}=\Bigl\{{1/2\atop 3/2}\ ,
\]
 $i=3:$
\[
\prod_{j=1}^{2}\left(1\mp\frac{1}{p_{j}}\right)=\left(1\mp\frac{1}{p_{1}}\right)\left(1\mp\frac{1}{p_{2}}\right)=\sum_{j=0}^{i-1}\sum_{\forall k}\frac{\left(\pm1\right)^{j}}{C_{j,i-1,k}}=\Bigl\{{1/3\atop 2}\ .
\]

Induktionsvoraussetzung: Die Formel gilt $\forall i\leq i_{0}-1=2$.
Dann gilt für $i=i_{0}:$ 

\begin{eqnarray*}
\prod_{j=1}^{i_{0}-1}\left(1\mp\frac{1}{p_{j}}\right) & = & \left(1\mp\frac{1}{p_{i_{0}-1}}\right)\prod_{j=1}^{i_{0}-2}\left(1\mp\frac{1}{p_{j}}\right)\\
 & = & \left(1\mp\frac{1}{p_{i_{0}-1}}\right)\sum_{j=0}^{i_{0}-2}\sum_{\forall k}\frac{\left(\pm1\right)^{j}}{C_{j,i_{0}-2,k}}=\sum_{j=0}^{i_{0}-1}\sum_{\forall k}\frac{\left(\pm1\right)^{j}}{C_{j,i_{0}-1,k}};
\end{eqnarray*}
\textbf{q.e.d.}

\subsection{Die Anzahl echt teilbarer natürlicher Zahlen $\mathbf{\mathbf{\leq x}}$ }

Die Zahl 1 gehört als multiplikative Invariante nicht zu den Primzahlen,
denn sie kann in beliebiger ganzzahliger Potenz einer beliebigen Zahl
hinzugefügt werden, ohne deren Wert zu ändern, obwohl sie der landläufigen
Primzahldefinition genügt, daß sie ganzzahlig nur durch 1 und sich
selbst geteilt werden kann. Sie soll im Folgenden mit den echt teilbaren
natürlichen Zahlen zusammen die Klasse der Nichtprimzahlen bilden.
Be\-zeich\-net $\sigma\left(x\right)$ die Anzahl aller Nichtprimzahlen
$\leq x$ und $\pi\left(x\right)$ die Anzahl aller Primzahlen $\leq x$
, so gilt
\begin{equation}
x=\sigma\left(x\right)+\pi\left(x\right).\label{eq:5.17}
\end{equation}
 \textbf{Satz 5.1.} Die Anzahl der echt teilbaren natürlichen Zahlen
(inclusive der Zahl 1) im abgeschlossenen Intervall $\left[1,x\right]$
beträgt mit $i_{0}$ gemäß Definition 5.4. und der eckigen Klammer
nach Definition 5.2. 
\begin{equation}
\sigma\left(x\right)=1+\sum_{i=1}^{i_{0}}\left(-1+\sum_{j=1}^{i-1}\left(-1\right)^{j}\sum_{\forall k}\left[\frac{x}{p_{i}C_{j,i-1,k}}\right]\right).\label{eq:5.18}
\end{equation}
 \textbf{Bew.:} Wie in allen Siebmethoden sollen in der Folge der
natürlichen Zahlen $\left\{ n_{i}\right\} ,\ n_{i}\leq x,$ sukzessive
alle Vielfachen der Primzahlen $p_{i}$ gestrichen und die Anzahlen
$\sigma_{i}\left(x\right)$ der echt durch $p_{i}$ teilbaren ermittelt
werden, die, beginnend mit der kleinsten Primzahl $p_{1}=2$, noch
ungestrichen stehen geblieben sind. 

Schritt 0: Streichung der natürlichen Zahl 1 in $\left\{ n_{i}\right\} $
als Nichtprimzahl. Es ist $\sigma_{0}\left(x\right)=1.$

Schritt 1: In $\left\{ n_{i}\right\} $ stehen $\sigma_{1}^{'}\left(x\right)=\left[\frac{x}{2}\right]$
Zahlen , die ohne Rest durch $p_{1}=2$ teilbar sind. Die erste dieser
Zahlen ist $p_{1}$ selbst, welche nicht zu streichen ist, sodaß $\sigma_{1}\left(x\right)=\left[\frac{x}{p_{1}}\right]-1$
Streichungen hinzukommen. Zusammen sind nun $s_{1}\left(x\right)=\sigma_{0}\left(x\right)+\sigma_{1}\left(x\right)=\left[\frac{x}{p_{1}}\right]$
Zahlen gestrichen. 

Schritt 2: In $\left\{ n_{i}\right\} $ waren $\sigma_{2}^{'}\left(x\right)=\left[\frac{x}{p_{2}}\right]$
Zahlen ohne Rest durch $p_{2}=3$ teilbar. Davon ist $p_{2}$ selbst
abzuziehen sowie alle geraden Vielfachen von $p_{2}$, die schon im
ersten Schritt erfaßt wurden, nämlich $\left[\frac{x}{p_{1}p_{2}}\right]$
Stück. Somit sind neu zu streichen $\sigma_{2}\left(x\right)=\left[\frac{x}{p_{2}}\right]-1-\left[\frac{x}{p_{1}p_{2}}\right]$,
sodaß gilt 
\[
s_{2}\left(x\right)=\sum_{i=0}^{2}\sigma_{i}\left(x\right)=1+\left[\frac{x}{p_{1}}\right]-1+\left[\frac{x}{p_{2}}\right]-1-\left[\frac{x}{p_{1}p_{2}}\right].
\]
 Unter Beachtung der Definitionen 5.2. und 5.3. erkennt man, daß die
Behauptung (\ref{eq:5.18}) bis zum $\left(i_{0}-1\right)-$ ten Schritt,
$\left(i_{0}-1\right)=2,$ erfüllt ist. Damit folgt durch Schluß von
$\left(i_{0}-1\right)$ auf $i_{0}$ 

Schritt $i_{0}:$ Die kleinste in den ersten $\left(i_{0}-1\right)$
Schritten nicht gestrichene Zahl, die noch $>p_{i_{0}-1}$ ist, ist
die nächste Primzahl $p_{i_{0}}$. Die Anzahl der in $\left\{ n_{i}\right\} $
existenten ganzzahligen Vielfachen von $p_{i_{0}}$ ist $\sigma_{i_{0}}^{'}\left(x\right)=\left[\frac{x}{p_{i_{0}}}\right]-1.$
Diese Zahl ist zu vermindern um alle bereits gestrichenen gemeinsamen
Vielfachen von $p_{i}$ mit $p_{j}\ ,j=1,2,\ldots,\left(i-1\right),$
also $\sum_{j=1}^{i-1}\left[\frac{x}{p_{i}p_{j}}\right].$ Diese letzte
Summe ist ihrerseits zu vermindern um die Vielfachen mit 3 gemeinsamen
Primfaktoren, also um $\sum_{j_{1,2}=1}^{i_{0}-1}\left[\frac{x}{p_{i_{0}}p_{j_{1}}p_{j_{2}}}\right],\ j_{2}>j_{1},$
die sonst doppelt gezählt würden. Allgemein ist jede solche Summe
mit k Faktoren im Nenner ihrerseits zu vermindern um eine solche mit
$\left(k+1\right)$ Nennerfaktoren, solange bis $k=i_{0}$ erreicht
wird:
\[
\sigma_{i_{0}}\left(x\right)=1+\sum_{i=1}^{i_{0}}\left(-1+\sum_{j=0}^{i-1}\left(-1\right)^{j}\sum_{\forall k}\left[\frac{x}{p_{i}C_{j,i-1,k}}\right]\right).
\]
Dies ist die behauptete Beziehung (\ref{eq:5.18}). Es bleibt lediglich
noch darauf hinzuweisen, daß das Verfahren mit dem $i_{0}-$ten Schritt
abbricht, wenn $i_{0}$ gemäß Definition 5.4. ermittelt wird, weil
spätestens dann alle Nichtprimzahlen $\leq x$ gestrichen sind;\textbf{
q.e.d.}

Zur Veranschaulichung der Gleichung (\ref{eq:5.18}) mögen folgende
Beispiele dienen.

\textbf{Beispiel 1: }$x=122,\ p_{i_{0}}=\prec\sqrt{122}\succ=11,$
also $i_{0}=5.$ Somit ergibt sich 

\begin{eqnarray*}
\sigma\left(x\right) & = & 1+\left(\left[\frac{x}{p_{1}}\right]-1\right)+\left(\left[\frac{x}{p_{2}}\right]-1-\left[\frac{x}{p_{1}p_{2}}\right]\right)
\end{eqnarray*}

\[
+\left(\left[\frac{x}{p_{3}}\right]-1-\left[\frac{x}{p_{1}p_{3}}\right]-\left[\frac{x}{p_{2}p_{3}}\right]+\left[\frac{x}{p_{1}p_{2}p_{3}}\right]\right)+s_{4}+s_{5}
\]

\begin{eqnarray*}
s_{4} & = & \left[\frac{x}{p_{4}}\right]-1-\left[\frac{x}{p_{1}p_{4}}\right]-\left[\frac{x}{p_{2}p_{4}}\right]-\left[\frac{x}{p_{3}p_{4}}\right]+\left[\frac{x}{p_{1}p_{2}p_{4}}\right]+\left[\frac{x}{p_{1}p_{3}p_{4}}\right]+\left[\frac{x}{p_{2}p_{3}p_{4}}\right]\\
 &  & -\left[\frac{x}{p_{1}p_{2}p_{3}p_{4}}\right]
\end{eqnarray*}

\begin{eqnarray*}
s_{5} & = & \left[\frac{x}{p_{5}}\right]-1-\left[\frac{x}{p_{1}p_{5}}\right]-\left[\frac{x}{p_{2}p_{5}}\right]-\left[\frac{x}{p_{3}p_{5}}\right]-\left[\frac{x}{p_{4}p_{5}}\right]+\left[\frac{x}{p_{1}p_{2}p_{5}}\right]+\left[\frac{x}{p_{1}p_{3}p_{5}}\right]
\end{eqnarray*}

\[
+\left[\frac{x}{p_{1}p_{4}p_{5}}\right]+\left[\frac{x}{p_{2}p_{3}p_{5}}\right]+\left[\frac{x}{p_{2}p_{4}p_{5}}\right]+\left[\frac{x}{p_{3}p_{4}p_{5}}\right]-\left[\frac{x}{p_{1}p_{2}p_{3}p_{5}}\right]
\]

\[
-\left[\frac{x}{p_{1}p_{2}p_{4}p_{5}}\right]-\left[\frac{x}{p_{1}p_{3}p_{4}p_{5}}\right]-\left[\frac{x}{p_{2}p_{3}p_{4}p_{5}}\right]+\left[\frac{x}{p_{1}p_{2}p_{3}p_{4}p_{5}}\right]\ ,
\]

\begin{eqnarray*}
\sigma\left(122\right) & = & 1+\left(61-1\right)+\left(40-1-20\right)+(24-1-12-8+4)\\
 &  & +\left(17-1-8-5-3+2+1+1-0\right)\\
 &  & +\left(11-1-5-3-2-1+1+1+0+\cdots+0\right)=92.
\end{eqnarray*}
 Daraus folgt nach (\ref{eq:5.17}) für die Anzahl Primzahlen bis
122 in Übereinstimmung mit Primzahltafeln $\pi\left(122\right)=122-\sigma\left(122\right)=30.$

Das Beispiel zeigt, daß es zur Berechnung von $\pi\left(122\right)$
reicht, alle Primzahlen\\
 $\leq\prec\sqrt{122}\succ=11=p_{5}$ zu kennen, obwohl die größte
Primzahl, die bis $x=122$ auftritt, $p_{30}=113$ lautet.

\textbf{Beispiel 2:} $x=168$ ; mittels derselben Formel erhält man
$\sigma\left(168\right)=129$, also $\pi\left(168\right)=39.$ Es
ist $p_{39}=167,$ denn erst ab $x=169$ ist $p_{6}=13$ in der Formel
für $\sigma\left(x\right)$ zu berücksichtigen. Für konkrete Berechnungen
von $\sigma\left(x\right)$ bzw. $\pi\left(x\right)$ für große Zahlen
x ist dieses Verfahren natürlich zu aufwändig. Ein möglichst frühzeitiges
Erkennen aller Nullsummanden wäre im Interesse einer Aufwandsreduzierung
wünschenswert, wie es im Falle des letzten Summanden möglich ist,
weil $\forall x>2$ gilt $p_{i}\downarrow>x.$ 

Offensichtlich bietet aber die Kenntnis der Gleichung (\ref{eq:5.18})
prinzipiell eine Möglichkeit, über die Primeigenschaft einer ungeraden
Zahl $x=2n+1$ zu entscheiden, indem man $\Delta=\sigma\left(2n+1\right)-\sigma\left(2n\right)$
bildet. Voraussetzung dafür ist die Kenntnis aller Primzahlen $\leq p_{i_{0}}=\prec\sqrt{2n+1}\succ.$
Ergibt sich $\Delta=1$, so ist x teilbar; im Falle $\Delta=0$ ist
x Primzahl. Es ist hierzu keinerlei Kenntnis über Primzahlen $>p_{i_{0}}$
erforderlich. Berechnet man zusätzlich $\sigma\left(2n+3\right)$,
so erfährt man, ob ein Primzahlzwilling vorliegt. Gilt $\Delta_{2}=\sigma\left(2n+3\right)-\sigma\left(2n+1\right)=2,$
so muß $x=2n+3$ neben $2n+2$ teilbar sein, für $\Delta_{2}=1$ ist
$x_{2}=2n+3$ Primzahl. D.h. für $\Delta=\Delta_{2}=1$ gilt $\left\{ 2n+1,2n+3\right\} =Primzahlzwilling.$
Auch hierzu zwei simple Beispiele aus dem Intervall $p_{3}^{2}=25<x<p_{4}^{2}=49.$
Wegen $i_{0}=3$ genügt die Betrachtung der Formel aus Beispiel 1
ohne $s_{4}$ und $s_{5}$.

\textbf{Beispiel 3:} $2n=28$ liefert $\sigma\left(2n\right)=19,\ \pi\left(2n\right)=9,\ \sigma\left(2n+1\right)=19,\ \pi\left(2n+1\right)=10.$
Also ist 29 eine Primzahl. Wegen $\sigma\left(2n+3\right)=20,\ \pi\left(2n+3\right)=11$
ist auch 31 Primzahl, also ist $\left\{ 29,31\right\} $ ein Primzahlzwilling.

\textbf{Beispiel 4:} $2n=40$ liefert $\sigma\left(40\right)=28,\ \pi\left(40\right)=12,$
$\sigma\left(41\right)=28,\ \pi\left(41\right)=13,\ \sigma\left(43\right)=29,\ \pi\left(43\right)=14,$
also ist auch $\left\{ 41,43\right\} $ Zwilling.

\subsection{Ein Satz zur Primzahlverteilung}

\textbf{Satz 5.2.} Sei n eine beliebige natürliche Zahl. Dann gibt
es $\forall n<\infty$ kein links offenes Intervall 
\begin{equation}
\eta_{n}:=\left(n^{2},\left(n+1\right)^{2}\right],\label{eq:5.19}
\end{equation}
 das nicht mindestens 2 Primzahlen enthält.

Es gilt sogar für die Anzahl der Primzahlen in $\eta_{n}$ 
\begin{equation}
\lim_{\eta_{n}\rightarrow\infty}\pi\left(\eta_{n}\right)\simeq\lim_{\eta_{n}\rightarrow\infty}\frac{2n+1}{2\ln\left(n+1\right)}=\infty.\label{eq:5.20.}
\end{equation}
 \textbf{Bew.: }Für $n=1$ enthält das Intervall $\eta_{1}=\left(1,4\right]$
nur 2 innere Zahlen, die Primzahlen 2 und 3. Die Intervallränder sind
konstruktionsbedingt stets nichtprim. Nach (\ref{eq:5.18}) ergibt
sich die Anzahl von Nichtprimzahlen im Intervall (\ref{eq:5.19})
zu 
\[
\sigma\left(\eta_{n}\right)=\sigma\left(\left(n+1\right)^{2}\right)-\sigma\left(n^{2}\right)=1+\sum_{i=1}^{i_{0}^{'}}\left(-1+\sum_{j=0}^{i-1}\left(-1\right)^{j}\sum_{\forall k}\left[\frac{\left(n+1\right)^{2}}{p_{i}C_{j,i-1,k}}\right]\right)
\]
 
\begin{equation}
-\left(1+\sum_{i=1}^{i_{0}}\left(-1+\sum_{j=0}^{i-1}\left(-1\right)^{j}\sum_{\forall k}\left[\frac{n^{2}}{p_{i}C_{j,i-1,k}}\right]\right)\right).\label{eq:5.21.}
\end{equation}
 Darin bedeuten $i_{0},i_{0}^{'}$ die Indizes aus 
\begin{equation}
p_{i_{0}}=\prec\sqrt{n^{2}}\succ=\prec n\succ\ \mbox{und}\ p_{i_{0}^{'}}=\prec n+1\succ.\label{eq:5.22}
\end{equation}
 Aus (\ref{eq:5.21.}) folgt 
\[
\sigma\left(\eta_{n}\right)=\sum_{i=1}^{i_{0}}\sum_{j=0}^{i-1}\left(-1\right)^{j}\sum_{\forall k}\left[\frac{\left(n+1\right)^{2}}{p_{i}C_{j,i-1,k}}\right]-\sum_{i=1}^{i_{0}}\sum_{j=0}^{i-1}\left(-1\right)^{j}\sum_{\forall k}\left[\frac{n^{2}}{p_{i}C_{j,i-1,k}}\right]+D
\]
 mit 
\begin{equation}
D=\sum_{i=i_{0}+1}^{i'_{0}}\left(-1+\sum_{j=0}^{i-1}\left(-1\right)^{j}\sum_{\forall k}\left[\frac{\left(n+1\right)^{2}}{p_{i}C_{j,i-1,k}}\right]\right).\label{eq:5.23}
\end{equation}
 Es ist 
\[
D=\Biggr\{{0\ \mbox{für}\ i_{0}=i_{0}^{'},\ \mbox{d.h.}\ n+1\neq Primzahl,\atop -1+\sum_{j=0}^{i_{0}}\left(-1\right)^{j}\sum_{\forall k}\left[\frac{\left(n+1\right)^{2}}{p_{i}C_{j,i-1,k}}\right]\ \mbox{sonst.}}
\]
 Da $n$ und $\left(n+1\right)$ benachbarte natürliche Zahlen sind,
kann $i_{0}^{'}>i_{0}$ nur gelten, wenn $\left(n+1\right)=p_{i_{0}+1}=p_{i'_{0}}$
selbst Primzahl ist. Für den Fall, daß $\left(n+1\right)$ nichtprim
ist, muß daher gelten $i_{0}^{'}=i_{0},$ woraus $D=0$ folgt. Dann
lassen sich die Summen vollständig zusammenfassen: 
\begin{equation}
\sigma\left(\eta_{n}\right)=\sum_{i=1}^{i_{0}}\sum_{j=0}^{i-1}\left(-1\right)^{j}\sum_{\forall k}\left(\left[\frac{\left(n+1\right)^{2}}{p_{i}C_{j,i-1,k}}\right]-\left[\frac{n^{2}}{p_{i}C_{j,i-1,k}}\right]\right).\label{eq:5.24.}
\end{equation}
 Am Schluß des Beweises von Satz 5.1. wurde darauf hingewiesen, daß
im Ausdruck für $\sigma\left(x\right)$ die Summe über $i$ automatisch
bei $i_{0}$ abbricht. Läßt man nun im Ausdruck $\sigma\left(\eta_{n}\right)$
auch für $\sigma\left(n^{2}\right)$ die Summe über $i$ bis $i_{0}^{'}$
laufen wie bei $\sigma\left(\left(n+1\right)^{2}\right),$ so begeht
man keinen Fehler, weil dadurch der Wert weder für $\sigma\left(n^{2}\right)$
noch für $\sigma\left(\eta_{n}\right)$ geändert wird. Ersetzt man
in (\ref{eq:5.24.}) $i_{0}$ durch $i_{0}^{'},$ so stellt diese
Gleichung den allgemein gültigen Ausdruck für $\sigma\left(\eta_{n}\right)$
dar. Wir lassen daher im Folgenden einfach D weg und interpretieren
das $i_{0}$ als $i_{0}^{'}.$ Durch Einführung der Summationsschrittweite
2 (in Zeichen: $\sum_{j=0(2)}^{i-1}\cdots$) lassen sich Terme mit
gleichen Vorzeichen in endlichen Summen zusammenfassen:
\begin{eqnarray}
\sigma\left(\eta_{n}\right) & = & \sum_{i=1}^{i_{0}}\Biggl\{\sum_{j=0\left(2\right)}^{i-1}\sum_{\forall k}\left[\frac{\left(n+1\right)^{2}}{p_{i}C_{j,i-1,k}}\right]-\sum_{j=1\left(2\right)}^{i-1}\sum_{\forall k}\left[\frac{\left(n+1\right)^{2}}{p_{i}C_{j,i-1,k}}\right]\nonumber \\
 &  & -\sum_{j=0\left(2\right)}^{i-1}\sum_{\forall k}\left[\frac{n^{2}}{p_{i}C_{j,i-1,k}}\right]+\sum_{j=1\left(2\right)}^{i-1}\sum_{\forall k}\left[\frac{n^{2}}{p_{i}C_{j,i-1,k}}\right]\Biggr\}.\label{eq:5.25}
\end{eqnarray}

Diese Gleichung sollte nun nach beiden Seiten abgeschätzt werden.
Mit Hilfe von (\ref{eq:5.9}) schätzt man ab $\left[\frac{\left(n+1\right)^{2}}{pC}\right]-\left[\frac{n^{2}}{pC}\right]\geq\left[\frac{2n+1}{pC}\right]$
und $\left[\frac{\left(n+1\right)^{2}}{pC}\right]-\left[\frac{n^{2}}{pC}\right]\leq-\left[-\frac{2n+1}{pC}\right],$
letzteres folgt durch Multiplikation mit $\left(-1\right)$ aus $\left[\frac{n^{2}}{pC}\right]-\left[\frac{\left(n+1\right)^{2}}{pC}\right]\geq\left[-\frac{2n+1}{pC}\right].$
Die Zusammenfassung beider Ungleichungen ergibt 
\begin{equation}
\left[\frac{2n+1}{pC}\right]\leq\left[\frac{\left(n+1\right)^{2}}{pC}\right]-\left[\frac{n^{2}}{pC}\right]\leq-\left[-\frac{2n+1}{pC}\right]\label{eq:5.26}
\end{equation}

und 
\begin{equation}
-\left[\frac{2n+1}{pC}\right]\geq-\left[\frac{\left(n+1\right)^{2}}{pC}\right]+\left[\frac{n^{2}}{pC}\right]\geq\left[-\frac{2n+1}{pC}\right].\label{eq:5.27}
\end{equation}
 In (\ref{eq:5.25}) werden nun die Terme mit gleichem Summationsschritt
zusammengefaßt sowie der 1. und 3. Term mittels (\ref{eq:5.26}),
der 2. und 4. Term mittels (\ref{eq:5.27}) abgeschätzt:
\begin{eqnarray}
\sigma\left(\eta_{n}\right) & \geq & \sum_{i=1}^{i_{0}}\left(\sum_{j=0\left(2\right)}^{i-1}\sum_{\forall k}\left[\frac{2n+1}{p_{i}C_{j,i-1,k}}\right]+\sum_{j=1\left(2\right)}^{i-1}\sum_{\forall k}\left[-\frac{2n+1}{p_{i}C_{j,i-1,k}}\right]\right)\nonumber \\
 & = & \sum_{i=1}^{i_{0}}\sum_{j=0}^{i-1}\sum_{\forall k}\left[\left(-1\right)^{j}\frac{2n+1}{p_{i}C_{j,i-1,k}}\right]\label{eq:5.28}
\end{eqnarray}
 sowie 
\begin{eqnarray}
\sigma\left(\eta_{n}\right) & \leq & \sum_{i=1}^{i_{0}}\left(\sum_{j=0\left(2\right)}^{i-1}\sum_{\forall k}-\left[-\frac{2n+1}{p_{i}C_{j,i-1,k}}\right]-\sum_{j=1\left(2\right)}^{i-1}\sum_{\forall k}\left[\frac{2n+1}{p_{i}C_{j,i-1,k}}\right]\right)\nonumber \\
 & =- & \sum_{i=1}^{i_{0}}\sum_{j=0}^{i-1}\sum_{\forall k}\left[\left(-1\right)^{j+1}\frac{2n+1}{p_{i}C_{j,i-1,k}}\right].\label{eq:5.29}
\end{eqnarray}
 Aus der Vereinigung von (\ref{eq:5.28}) und (\ref{eq:5.29}) folgt
nach Abschätzung gemäß (\ref{eq:5.8.}) - imFalle von (\ref{eq:5.29})
mit umgekehrtem Vorzeichen - 
\begin{equation}
\sum_{i=1}^{i_{0}}\sum_{j=0}^{i-1}\sum_{\forall k}\left(-1\right)^{j}\frac{2n+1}{p_{i}C_{j,i-1,k}}\leq\sigma\left(\eta_{n}\right)\leq\sum_{i=1}^{i_{0}}\sum_{j=0}^{i-1}\sum_{\forall k}\left(-1\right)^{j}\frac{2n+1}{p_{i}C_{j,i-1,k}}\ ,\label{eq:5.30}
\end{equation}
 also die Gleichheit. Durch Einsetzen von (\ref{eq:5.15}) ergibt
sich 
\begin{equation}
\sigma\left(\eta_{n}\right)=\left(2n+1\right)\sum_{i=1}^{i^{0}}p_{i}^{-1}\sum_{j=0}^{i-1}\sum_{\forall k}\frac{\left(-1\right)^{j}}{C_{j,i-1,k}}=\left(2n+1\right)\sum_{i=1}^{i_{0}}p_{i}^{-1}\prod_{j=1}^{i-1}\left(1-\frac{1}{p_{j}}\right).\label{eq:5.31}
\end{equation}
 Darin gilt nach Definition $\prod_{j=1}^{0}\left(1-\frac{1}{p_{j}}\right)=1.$
Wegen (\ref{eq:5.17}), angewendet auf das Intervall $\eta_{n}$ mit
der Breite $b_{n}=\left(n+1\right)^{2}-n^{2}=2n+1,$ gilt dann 
\begin{equation}
\pi\left(\eta_{n}\right)=2n+1-\sigma\left(\eta_{n}\right)=\left(2n+1\right)\left(1-\sum_{i=1}^{i_{0}}p_{i}^{-1}\prod_{j=1}^{i-1}\left(1-\frac{1}{p_{j}}\right)\right).\label{eq:5.32.}
\end{equation}
 Ersetzt man $-p_{i}^{-1}=\left(-1+\left(1-p_{i}^{-1}\right)\right)$
und multipliziert dies aus (Umordnungen in endlichen Produkten sind
erlaubt) , so heben sich Summanden paarweise auf. Man erhält 
\begin{equation}
S_{n}:=\frac{\pi\left(\eta_{n}\right)}{2n+1}=1+\sum_{i=1}^{i_{0}}\left(\prod_{j=1}^{i}\left(1-\frac{1}{p_{j}}\right)-\prod_{j=1}^{i-1}\left(1-\frac{1}{p_{j}}\right)\right)=\prod_{j=1}^{i_{0}}\left(1-\frac{1}{p_{j}}\right).\label{eq:5.33.}
\end{equation}
 Diese Formel liefert bei numerischer Auswertung für kleine Werte
von $n,\ p_{i_{0}}=\prec n+1\succ$ die reale Anzahl von Primzahlen
im n-ten Quadratintervall mit einem maximalen Fehler von $\pm2$ für
$n\leq40$ bei Verwendung ganzzahlig abgerundeter Werte. Auch die
daraus bestimmte Anzahl von Primzahlen $\leq x$, $\pi\left(x\right)=\sum_{n=0}^{n_{0}}\pi\left(\eta_{n}\right),$
liefert, s. \textbf{Tab. 5.1.}, sinnvolle Werte: 

~~\ \quad{}%
\begin{tabular}{|c|c|c|c|c|}
\hline 
$n_{0}$ & 10 & 20 & 30 & 40\tabularnewline
\hline 
$\pi\left(\left(n_{0}+1\right)^{2}\right)$ & 26 & 80 & 161 & 266\tabularnewline
\hline 
$Realwert$ & 30 & 85 & 162 & 263\tabularnewline
\hline 
\end{tabular}

Wir bilden nun $\ln S_{n}$, spalten den 1. Term ab und entwickeln
dann den Logarithmus in eine Potenzreihe: 
\begin{eqnarray}
\ln S_{n} & = & \ln\prod_{j=1}^{i_{0}}\left(1-\frac{1}{p_{j}}\right)=\ln\frac{1}{2}+\sum_{j=2}^{i_{0}}\ln\left(1-\frac{1}{p_{j}}\right)=\ln\frac{1}{2}-\sum_{j=2}^{i_{0}}\sum_{\nu=1}^{\infty}\frac{1}{\nu p_{j}^{\nu}}\nonumber \\
 & = & \ln\frac{1}{2}-\sum_{j=2}^{i_{0}}\frac{1}{p_{j}}-\sum_{j=2}^{i_{0}}\sum_{\nu=2}^{\infty}\frac{1}{\nu p_{j}^{\nu}}.\label{eq:5.34}
\end{eqnarray}
 Darin ist $A=-\sum_{j=2}^{i_{0}}\sum_{\nu=2}^{\infty}\frac{1}{\nu p_{j}^{\nu}}=const.$,
denn wegen $\nu\geq2$ ist A eine absolut konvergente Reihe. Aus der
Literatur \cite{key-1}, p.343, Satz 5, entnehmen wir die Aussage
\begin{equation}
\sum_{j=1}^{i_{0}}\frac{1}{p_{j}}=\ln\ln p_{i_{0}}+c+O\left(\frac{1}{\ln p_{i_{0}}}\right),\label{eq:5.35}
\end{equation}
 worin $c=const.$ und $O\left(\frac{1}{\ln p_{i_{0}}}\right)$ den
Fehlerterm bezeichnet. Der Wert $\ln\ln x=0$ wird für $x=e$ angenommen,
$\ln\ln p_{1}=\ln\ln2=-0,3665129\ldots<0,\ \ln\ln3=0,094047827\ldots,$
$\ln3=1,098612289\ldots$ . Da in (\ref{eq:5.34}) $j=2,\ldots,i_{0}$
läuft, haben wir die Belegung der Konstanten etwas zu modifizieren.
Um den Wert von $\ln\ln p_{i_{0}}$ benutzen zu können, haben wir
davon $\ln\ln3$ abzuziehen. So erhält man schließlich 
\begin{equation}
\ln S_{n}=\ln\frac{1}{2}-\left(\ln\ln p_{i_{0}}-\ln\ln3+c+O\left(\frac{1}{\ln p_{i_{0}}}\right)+A\right)\label{eq:5.36}
\end{equation}
sowie durch Exponentiation 
\begin{equation}
S_{n}=\frac{\pi\left(\eta_{n}\right)}{2n+1}=\frac{1}{2}\cdot\frac{1}{\ln p_{i_{0}}}\cdot A^{'}\label{eq:5.37}
\end{equation}
mit $A'=\left(\ln3\right)e^{-c+A-O\left(1/\ln p_{i_{0}}\right)}.$
Dies ist die Aussage des Gauß'schen Primzahlsatzes für Quadratintervalle.
Für $0\simeq A-c-O\left(1/\ln p_{i_{0}}\right)$ wird $A^{'}\simeq\ln3.$ 

Wenn aber die Primzahlverteilung über jedem Quadratintervall diesem
logarithmi-schen Gesetz gehorcht, gilt sie für jedes beliebige größere
Intervall, also auch für $\left(1,x\right],$ wie es der Gauß'sche
Primzahlsatz verlangt: $\pi\left(x\right)=A^{''}\left(\frac{x}{\ln x}\right)\ln3,$
denn $2\ln p_{i_{0}}=\ln p_{i_{0}}^{2}\simeq\ln\left(n+1\right)^{2}=\ln x.$
Die Formel $\pi\left(x\right)=\left(\frac{x}{\ln x}\right)$ gibt
im allgemeinen Werte, die um $5\ldots10\%$ zu niedrig liegen. Der
Faktor $\ln3$ vergrößert sie um 9,86\% , sodaß wohl $A^{''}<1$,
aber nahe bei 1, anzunehmen ist.

Wir haben hier das logarithmische Verteilungsgesetz für Quadratintervalle
abgeleitet aus einer (ziemlich koplizierten) exakten Formel, daraus
auf seine allgemeine Gültigkeit geschlossen und den Gauß'schen Satz
als Approximation erhalten; \textbf{q.e.d.}

Die \textbf{Abb. 5.1.} zeigt $\forall n\leq240$ die Anzahlen $\pi\left(\eta_{n}\right)=$PiEn
und $\pi\left(\eta_{A,n}\right)=\pi\left(2n+1\right)=\frac{A^{'}\left(2n+1\right)}{\ln\left(n+1\right)}$
=PiEAn sowie die gemäß (\ref{eq:5.40.}) streuenden Funktionen $\pi_{s}\left(\eta_{n}\right)=\pi\left(\eta_{n}\right)\left(1+\delta\right)$
mit $\delta=0$ für den Mittelwert $g\left(n\right)$, $\delta=+\frac{A^{'}}{\ln\left(n+1\right)}$
als obere Schranke $go\left(n\right)$, $\delta=-\frac{A^{'}}{\ln\left(n+1\right)}$
als untere Schranke $gu\left(n\right)$. Der im Mittel stetige Anstieg
ist deutlich zu erkennen, ebenso daß im Anfangsintervall gleicher
Größe $\eta_{A,n}:=\left[0,2n+1\right]$ der Gauß'sche Satz $\pi\left(\eta_{A,n}\right)=\frac{A^{'}\left(2n+1\right)}{\ln\left(n+1\right)}\approx2\pi\left(\eta_{n}\right)$
ergibt.

\begin{landscape}
\begin{samepage}

\begin{figure}[th]
\caption{ 5.1.: Primzahlanzahlen über Quadratintervallen $\pi(\eta_{n})$ und
in gleichgroßen Anfangsintervallen $\pi(\eta_{A,n})$}

\includegraphics{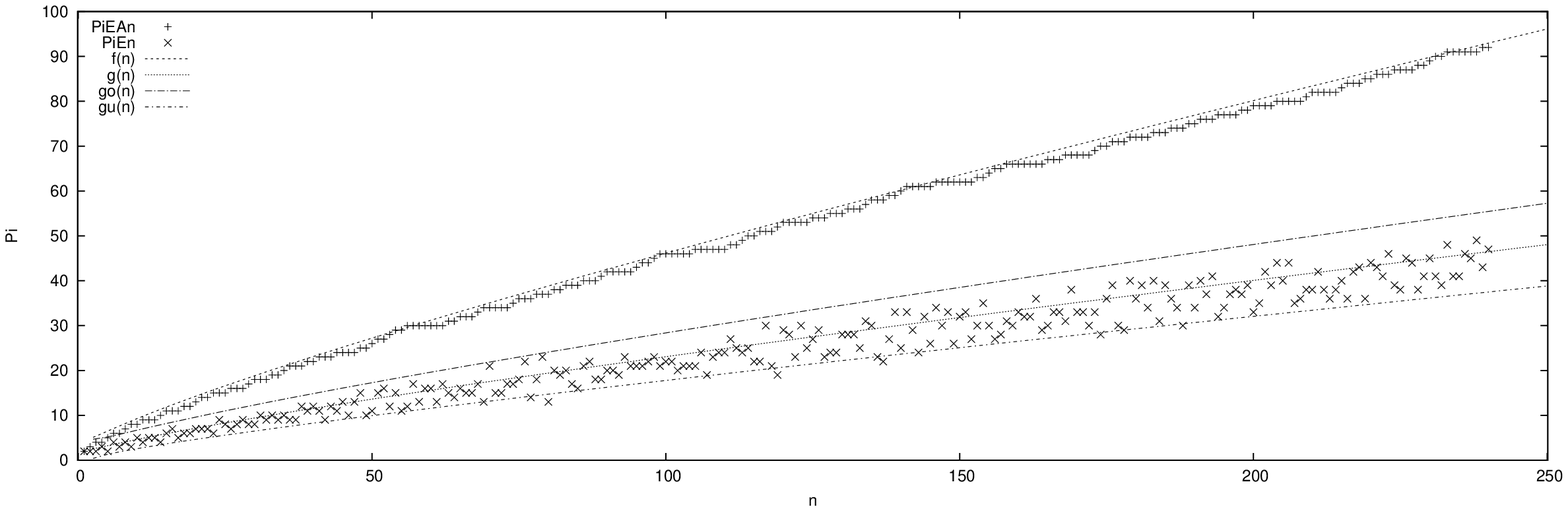}

\begin{minipage}[t]{1\columnwidth}%
\caption{9.1.: Streuende Anzahl von Goldbachpaar-Darstellungen $\pi_{g}(2m)$}

\includegraphics{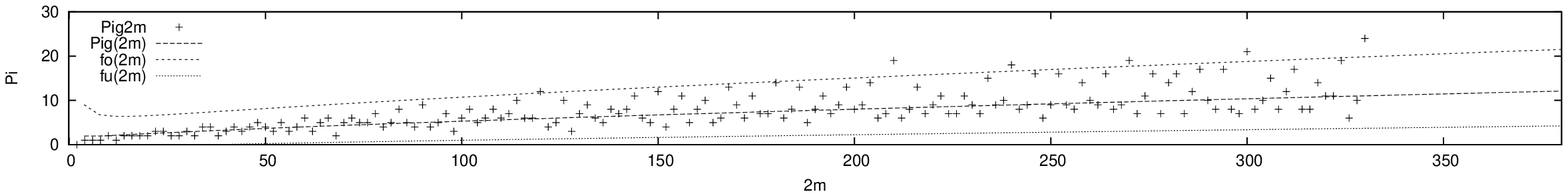}%
\end{minipage}
\end{figure}

\end{samepage}
\end{landscape}
\newpage

\subsection{Vertrauensgrenzen des Verteilungssatzes}

Nach (\ref{eq:5.37}) haben wir mit $p_{i_{0}}=\prec n+1\succ\cong n+1$
als Wahrscheinlichkeit W aus mittlerer Anzahl $\pi\left(\eta_{n}\right)$
der Primzahlen im Intervall $\eta_{n}$ mit der Breite $b_{n}=2n+1$
erhalten:
\begin{equation}
W=\frac{\pi\left(\eta_{n}\right)}{b_{n}}=\frac{A^{'}}{2\ln\left(n+1\right)}\approx\frac{A^{'}}{\ln x}\ \mbox{mit}\ x\epsilon\eta_{n}.\label{eq:5.38}
\end{equation}
 Der Kehrwert $L_{n}=\frac{1}{W}=\frac{b_{n}}{\pi\left(\eta_{n}\right)}$
kann als mittlerer Platzbedarf einer Primzahl in $\eta_{n}$ gedeutet
werden. Wir bestimmen stattdessen den halben Platzbedarf $\frac{1}{2}L_{n}^{'}$
eines Primzahlpaares. Die Wahrscheinlichkeit, daß 2 Zahlen Primzahlen
in $\eta_{n}$ sind, wird durch ihr Wahrscheinlichkeitsprodukt $W_{ges.}=W_{p_{1}}W_{p_{2}}=W^{2}$
gegeben. Dann ist die Streubreite für eine Primzahl in $\eta_{n}$
\begin{equation}
S_{n}=\frac{b_{n}}{\frac{1}{2}L_{n}^{'}}=2b_{n}W^{2}=2\left(2n+1\right)\left(\frac{A^{'}}{2\ln\left(n+1\right)}\right)^{2}.\label{eq:5.39.}
\end{equation}

Da $p_{1}=2$ die kleinste Primzahl ist, kann (außer zwischen $p_{1}$
und $p_{2}=3$) der minimale Primzahlabstand nicht kleiner als 2 sein.
Der Abstand zwischen zwei benachbarten Primzahlen ist kleiner als
das zugehörige $b_{n}$, deshalb gilt $2\leq p_{i+1}-p_{i}<b_{n}.$

Deshalb erhalten wir für die in $\eta_{n}$ streuende Anzahl von Primzahlen

\begin{equation}
\pi_{s}\left(\eta_{n}\right)=b_{n}W\pm S_{n}=\frac{A^{'}\left(n+\frac{1}{2}\right)}{\ln\left(n+1\right)}\pm\frac{A^{'2}\left(n+\frac{1}{2}\right)}{\ln^{2}\left(n+1\right)}=\frac{A^{'}\left(n+\frac{1}{2}\right)}{\ln\left(n+1\right)}\left(1\pm\frac{A^{'}}{\ln\left(n+1\right)}\right).\label{eq:5.40.}
\end{equation}

bzw. $\pi_{s}\left(x\right)=\frac{A^{'}x}{\ln x}\pm\frac{2A^{'2}x}{\ln^{2}x}\ .$
Dieser Toleranzbereich ist realistischer als die willkürliche Annahme
$\pi_{s}\left(x\right)=\frac{A^{'}x}{\ln x}\left(1\pm0,2\right)$
und beschreibt die Realität gut, s. Abb. 5.1. und \textbf{}\\
\textbf{Tab. 5.2.}:

\begin{tabular}{|c|c|c|c|c|c|c|c|c|c|c|}
\hline 
n & 1 & 10 & 100 & 150 & 300 & 400 & $10^{3}$ & $10^{4}$ & $10^{5}$ & $10^{6}$\tabularnewline
\hline 
$\pi\left(\eta_{n}\right)$ & 2,29 & 4,4 & 23,1 & 31,9 & 55,8 & 70,8 & 154 & 1151 & 9207 & 76725,4\tabularnewline
\hline 
$L_{n}$ & 0,86 & 10,2 & 37,9 & 44,8 & 58,0 & 64,0 & 85,0 & 151,0 & 235,9 & 339,7\tabularnewline
\hline 
$\pm\frac{2n+1}{L_{n}}$ & 3,5 & 2,05 & 5,30 & 6,72 & 10,4 & 12,5 & 23,6 & 132,5 & 847,7 & 5886,9\tabularnewline
\hline 
\end{tabular}

Größere Schwankungen sind möglich: Der wirkliche Extremfall, daß die
erste Primzahl am unteren, die zweite am oberen Intervallrand liegt,
könnte schlimmstenfalls das Doppelte vom mittleren Primzahlabstand
erbringen und statistische ``Ausreißer'' erklären. 

Bezeichnet $\Delta\pi_{max}\left(x\right)$ die in der weiteren Umgebung
von x auftretenden maximalen Primzahlabstände, so gibt \textbf{Tab
5.3.} den Vergleich mit den Primzahlabständen gemäß (\ref{eq:5.39.}
):

\begin{tabular}{|c|c|c|c|c|c|c|}
\hline 
$x$ & 10 & 10$^{2}$ & 10$^{3}$ & 10$^{4}$ & 10$^{5}$ & 8,4$\cdot10^{5}$\tabularnewline
\hline 
$\Delta\pi_{max}\left(x\right)$ & 4 & 14 & 20 & 36 & 54 & 100\tabularnewline
\hline 
$\frac{\ln^{2}x}{2A^{'2}}$ & 2,36 & 9,44 & 21,23 & 37,75 & 58,98 & 82,81\tabularnewline
\hline 
\end{tabular} 

Es wurde mit $A^{'}=1,06$ gerechnet.

\section{Ein falscher Satz}

In der Literatur findet man seit langer Zeit die Behauptung, man könne
``beliebig große primzahlfreie Intervalle'' konstruieren (s. z.B.
\cite{key-1}, p. 22 Zeile 12) und die Frage nach einer oberen Schranke
für den Abstand benachbarter Primzahlen sei sinnlos. Zum Beweis konstruiert
man die Intervalle $\left[n!\pm2,n!\pm n\right]$ , die primzahlfreie
Intervalle der Mindestlänge $\left(n-1\right)$ darstellen, und läßt
$n\rightarrow\infty$ laufen. Von jeder endlichen Zahl $n!-1$ rückwärts
gezählt, steht natürlich kein beliebig $\left(=\infty\right)$ großes
primzahlfreies Intervall zur Verfügung. Von $n!+1$ nach ``oben''
gezählt, steht zunächst nur das endliche primzahlfreie Intervall der
Länge $n-1$ bereit. ( De facto kann das primzahlfreie Intervall auch
doppelt so groß sein!) Die angegebene Konstruktion beweist in der
Tat nur, daß ``im unendlich fernen Punkt'' selbst genau ein beliebig
großes primzahlfreies Intervall gedacht werden kann. Dieser ``Punkt``
ist aber singulär und zugleich der einzige Häufungspunkt der Primzahlen.
Der Satz ist also falsch.

Wegen Satz 5.2. gibt es kein primzahlfreies Quadratintervall $\eta_{n}=\left(n^{2},\left(n+1\right)^{2}\right].$
Daraus folgt, daß der maximale Abstand zwischen benachbarten Primzahlen
$<2n$ für jede vorgelegte Zahl $x=n^{2}$ ist, d.h. $\Delta_{max}=p_{i+1}-p_{i}<2\left[\sqrt{p_{i+1}}\right].$
Man kann sogar den Mittelwert $\Delta_{max,M}$ sowie obere und untere
Schranken dazu aus dem logarithmi\-schen Verteilungsgesetz über Quadratintervallen
ableiten. Wir werden in den folgenden Kapiteln dieser Frage nochmals
begegnen. Solche Angaben werden natürlich n-abhängig sein. 

Beim Vergleich verschieden mächtiger Mengen muß etwas mehr Sorgfalt
aufgewendet werden. Eine Zahl $m=p_{n}\downarrow$ oder $m=n!$ gehört
zur Potenzordnung n, da sie n Faktoren enthält. Das Quadratintervall
$\left(\sqrt{m}^{2},\left(\sqrt{m}+1\right)^{2}\right]$, in dem m
liegt, hat die Potenzordnung $2<n$ und die lineare Breite $b=2\sqrt{m}+1=1+2\sqrt{n!}\gg n-1=$
primzahlfreies Intervall, welches eine lineare Mannigfaltigkeit darstellt,
das schon die umgebende Mannigfaltigkeit 2. Ordnung (Quadratintervall)
nicht auszufüllen vermag. Dies dürfte die Ursache des Trugschlusses
im angeführten Satz von der Konstruierbarkeit ``beliebig großer''
primzahlfreier Intervalle gewesen sein.

Obwohl also zu jeder beliebig gewählten natürlichen Zahl n<$\infty$
eine natürliche Zahl n'=n!+1konstruierbar ist, der ein primzahlfreies
Intervall der Mindestlänge n-1 folgt, gibt es keine natürliche Zahl,
der ein beliebig langes primzahlfreies Intervall folgen kann.

\section{Zur Verteilung von Primzahlzwillingen}

Die Primzahlen sind bezüglich der Multiplikation statistisch unabhängig,
denn eine jede Primzahl p enthält als Teiler nur 1 und p und ist selbst
kein Teiler irgend einer anderen Primzahl. Primzahlen sind sozusagen
die Ein-heiten der Produktmengen. Es ist daher legitim, die Zwillingsbedingung
$p_{i+1}=p_{i}+2$ zu kombinieren mit dem Primzahlverteilungsgesetz
über einem Quadratintervall $\eta_{n}$ oder einem größeren Intervall.
Zuvor stellen wir noch fest, daß ein Primzahlzwilling mit $p_{i}>3$
niemals eine Quadratzahl umgreifen kann. Es müssen daher - außer $\left\{ p_{2}=3,p_{3}=5\right\} $
mit $2^{2}$ in der Mitte - die beiden Partner eines beliebigen Primzahlzwillings
$>5$ stets vollständig ein und demselben Quadratintervall angehören. 

\textbf{Bew.:} $x^{2}-1=\left(x+1\right)\left(x-1\right)$ ist $\forall i>2$
eine stets teilbare Zahl, deshalb kann höchstens $x^{2}+1$ Primzahl
sein, sodaß zwar ein Primzahlzwilling am Anfang eines Quadratintervalls,
aber vollständig darin stehen kann, wenn x geradzahlig ist. Ist x
ungerade, so stehen 2 teilbare Zahlen nebeneinander: $x^{2}$und $\left(x-1\right)\left(x+1\right),$
sodaß ebenso nur der Anfang $x^{2}+1$ des größeren Quadratintervalls
für die Bildung eines Primzahlzwillings in Betracht kommt; \textbf{q.e.d.}

Für die Primzahldichte über einem Quadratintervall (wahrscheinlichster
oder Mittelwert) ergab sich in (\ref{eq:5.37}) 
\begin{equation}
W=\frac{\pi\left(\eta_{n}\right)}{2n+1}=\frac{A^{''}\ln3}{2\ln p_{i_{0}}}\label{eq:7.1.}
\end{equation}
 mit $p_{i_{0}}=\prec n+1\succ,\ \eta_{n}=\left(n^{2},\left(n+1\right)^{2}\right].$
Ein Primzahlzwilling ist gekennzeichnet durch die Annahme des geringstmöglichen
Primzahlabstandes 2, also $p_{i+1}=p_{i}+2.$ Da $p_{i},p_{i+1}$
im gleichen Quadratintervall liegen müssen, wird ihre Wahrscheinlichkeit
durch dasselbe $p_{i_{0}}$ bestimmt. Es darf also die Zwillingswahrscheinlichkeitsverteilung
als $W_{z}=W_{i}W_{i+1}=W^{2}$ angesetzt werden. Bezeichnet $\pi_{2}\left(\eta_{n}\right)$
die mittlere oder wahrscheinliche Anzahl von Primzahlzwillingen in
$\eta_{n},$ so kann geschrieben werden
\[
W_{z}=\frac{\pi_{2}\left(\eta_{n}\right)}{2n+1}=\frac{A^{''2}\ln^{2}3}{4\ln^{2}p_{i_{0}}}
\]
 oder
\begin{equation}
\pi_{2}\left(\eta_{n}\right)=\frac{\left(2n+1\right)A^{''2}\ln^{2}3}{4\ln^{2}p_{i_{0}}}\approx\frac{\left(n+\frac{1}{2}\right)A^{''2}\ln^{2}3}{2\ln^{2}\left(n+1\right)}.\label{eq:7.2.}
\end{equation}
 Diese Funktion ist in \textbf{Abb. 7.1.} zwischen die Punkte der
numerisch bestimmten Prim\-zahl\-zwillings-Anzahlen je $\eta_{n}\ \forall n\leq915$
eingezeichnet und zeigt eine gute Approximation an. Es wurde $A^{'}=A^{''}\ln\left(3\right)=1,06$
benutzt.

\textbf{Satz 7.1.} Die mittlere Anzahl von Primzahlzwillingen im Quadratintervall
$\eta_{n}$ steigt gemäß (\ref{eq:7.2.}) mit n monoton an und wächst
für $n\rightarrow\infty$ über alle Schranken. Die integrale Anzahl
der Primzahlzwillinge wächst dann natürlich erst recht über alle Schranken.

\textbf{Bew.:} Aus der Divergenz des Ausdrucks $y=\frac{\sqrt{x-1/2}}{\ln x}$
(der Logarithmus wächst schwächer als jede Potenzfunktion $x^{\nu}\ mit\ \nu>0$
) folgt mit $x=n+1$ auch die Divergenz von $y^{2}=\frac{n+1/2}{\ln^{2}\left(n+1\right)},$
daran ändert auch der konstante Faktor $\frac{\ln^{2}3}{2}$ nichts;
\textbf{q.e.d.}

Dies Ergebnis ist insofern erstaunlich, als im Anfangsbereich ($n\leq122$)
primzahlleere Quadratintervalle existieren für $n=\{9,19,26,27,30,34,39,49,53,77,122\}$,
darunter sogar 2 benachbarte Quadratintervalle. Fragt man aber nach
dem größten Wert von n, oberhalb dessen ein vorgegebener ganzzahliger
Wert von $\pi_{2}\left(\eta_{n}\right)$ nicht mehr unterschritten
wird, so ergibt sich mit $\pi_{2}\left(\eta_{n,max}\right)$ nach
(\ref{eq:7.2.}) aus der Abb.7.1. folgende Tabelle \textbf{Tab.7.1.}: 

~\quad{}\qquad{}%
\begin{tabular}{|c|c|c|c|c|c|}
\hline 
\selectlanguage{english}%
$n_{max}$\selectlanguage{ngerman}%
 & 122 & 213 & 502 & 545 & 829\tabularnewline
\hline 
$\pi_{2}\left(\eta_{n}\right)$ & 0 & 1 & 2 & 3 & 4\tabularnewline
\hline 
$\pi_{2}\left(\eta_{n,max}\right)$ & 3,19 & 4,47 & 7,84 & 8,29 & 11,09\tabularnewline
\hline 
\end{tabular}

\begin{landscape}

\psfrag{Pi2}{\begin{math}\pi_2\end{math}}

\begin{figure}[th]

\caption{7.1.: $\pi_{2}(\eta_{n})$ experimentell und theoretisch als $f(n)$
incl. Streubereich}

\includegraphics{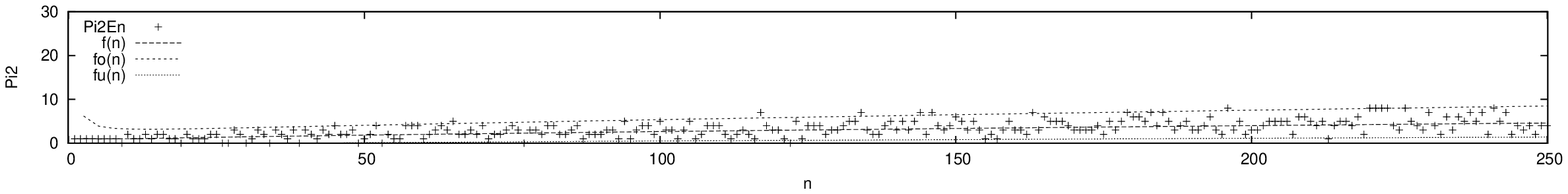}

\includegraphics{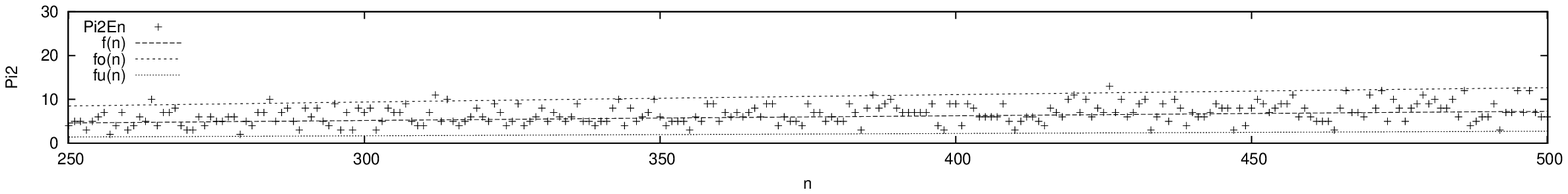}

\includegraphics{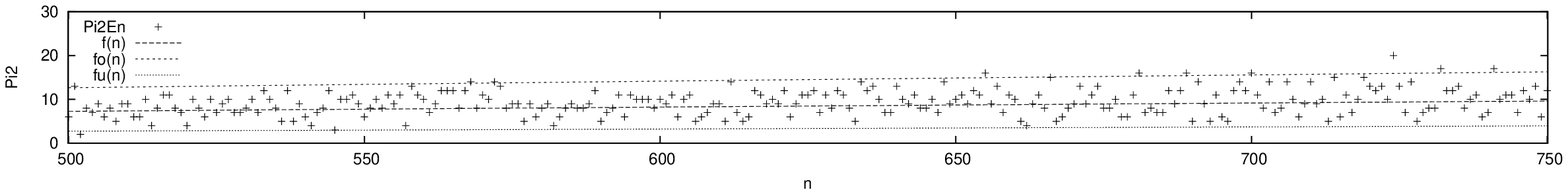}

\includegraphics{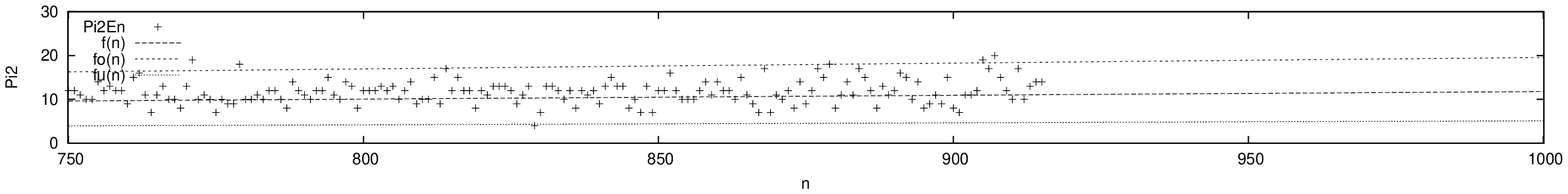}
\end{figure}

\end{landscape}
\newpage

Das wirft die Frage nach der Streubreite um die mittlere theoretische
Häufigkeit $\pi_{2}\left(\eta_{n}\right)$ auf. Benutzt man in (\ref{eq:5.37})
für $A^{'}=1,06$ anstelle von $A^{''}\ln3=A^{''}\cdot1,0986123\ldots$,
so stimmen der Mittelwert
\begin{equation}
\pi\left(\eta_{n}\right)=\frac{n+\frac{1}{2}}{\ln\left(n+1\right)}\cdot1,06\ ,\label{eq:7.3.}
\end{equation}
die obere Schranke 
\begin{equation}
\pi_{ob.}\left(\eta_{n}\right)=1,2\cdot1,06\cdot\frac{n+\frac{1}{2}}{\ln\left(n+1\right)}\label{eq:7.4.}
\end{equation}
 und die untere Schranke 
\begin{equation}
\pi_{unt.}\left(\eta_{n}\right)=0,8\cdot1,06\cdot\frac{n+\frac{1}{2}}{\ln\left(n+1\right)}\label{eq:7.5.}
\end{equation}
für $n<1000$ befriedigend mit der Realität überein. Durch Quadrieren
der Wahrscheinlichkeitsausdrücke ergibt sich damit für Primzahlzwillinge
\begin{equation}
\pi_{2}\left(\eta_{n}\right)=\frac{n+\frac{1}{2}}{2\ln^{2}\left(n+1\right)}\cdot1,06^{2}\ ,\label{eq:7.6.}
\end{equation}
\begin{equation}
\pi_{2,ob.}\left(\eta_{n}\right)=\frac{n+\frac{1}{2}}{2\ln^{2}\left(n+1\right)}\left(1,2\cdot1,06\right)^{2}\ ,\label{eq:7.7.}
\end{equation}
\begin{equation}
\pi_{2,unt.}\left(\eta_{n}\right)=\frac{n+\frac{1}{2}}{2\ln^{2}\left(n+1\right)}\left(0,8\cdot1,06\right)^{2}\ .\label{eq:7.8.}
\end{equation}
Während der Mittelwert gut liegt, ist die reale Abweichung vom Mittelwert
etwa doppelt so groß. Dieser Mangel wird überwunden, wenn man statt
der willkürlichen Streufaktoren 0,8 und 1,2 die besser begründete
mittlere Abweichung gemäß (\ref{eq:5.40.}) einführt. 

Wir betrachten 2 Primzahlzwillinge mit der Gesamtwahrscheinlichkeit
$W_{ges.}=W_{z_{1}}W_{z_{2}}=W_{s}^{2},$ weil für genügend große
n der Zwillingsabstand $<b_{n}$ ist. Für den einzelnen Primzahlzwilling
verdoppeln wir wieder den streuenden Anteil und erhalten 
\begin{eqnarray}
\pi_{2,s}\left(\eta_{n}\right) & = & \left(2n+1\right)\left(\frac{A^{'}}{2\ln\left(n+1\right)}\right)^{2}\left(1+2\left(4\left(\frac{A^{'}}{2\ln\left(n+1\right)}\right)^{2}\pm4\frac{A^{'}}{2\ln\left(n+1\right)}\right)\right)\nonumber \\
 & = & \left(2n+1\right)\left(\frac{A^{'}}{2\ln\left(n+1\right)}\right)^{2}\left(1+\delta_{2}\right),\label{eq:7.9.}
\end{eqnarray}

Darin sind die Streufaktoren $\delta_{2}=0$ für den Mittelwert, für
die obere und untere Schranke 
\begin{equation}
\delta_{2}=\delta_{2,+}=2\left(\frac{A^{'}}{\ln\left(n+1\right)}\right)^{2}+\frac{4A^{'}}{\ln\left(n+1\right)}\ ,\label{eq:7.10.}
\end{equation}
 
\begin{equation}
\delta_{2}=\delta_{2,-}=2\left(\frac{A^{'}}{\ln\left(n+1\right)}\right)^{2}-\frac{4A^{'}}{\ln\left(n+1\right)}\ .\label{eq:7.11.}
\end{equation}
 Die 3 Funktionen sind in Abb. 7.1. neben den Werten des numerischen
Experiments eingezeichnet; sie beschreiben die Realität hinreichend
gut. Die Tab. 7.2. zeigt, daß die nach (\ref{eq:7.9.}) berechnete
gesamte Streubreite sich oberhalb eines Schwellwertes über die n-Achse
erhebt und für große n nahezu symmetrisch wird, während bei kleinen
n deutlich $\delta_{2,+}>\delta_{2,-}$ gilt. \textbf{Tab. 7.2.}: 

\begin{tabular}{|c|c|c|c|c|c|c|c|}
\hline 
n & 1 & 10 & $10^{2}$ & $10^{3}$ & $10^{4}$ & $10^{5}$ & $10^{6}$\tabularnewline
\hline 
$\pi_{2}\left(\eta_{n}\right)$ & 1,75 & 1,03 & 2,65 & 11,78 & 66,23 & 423,85 & 2943,39\tabularnewline
\hline 
$\pi_{2}\left(\eta_{n}\right)\delta_{2,+}$ & 11,19 & 2,14 & 4,01 & 15,68 & 82,32 & 505,23 & 3411,39\tabularnewline
\hline 
$\pi_{2}\left(\eta_{n}\right)\delta_{2,-}$ & 0,49 & 0,32 & 1,57 & 8,45 & 51,86 & 349,25 & 2510,71\tabularnewline
\hline 
\end{tabular}

\section{Primzahlvierlinge}

\begin{wrapfigure}{o}{0.5\columnwidth}%
\centering

\caption{8.1. $\pi_{4}(n^{4})$ und $\pi_{4}(\eta_{4,n})$ }

\includegraphics{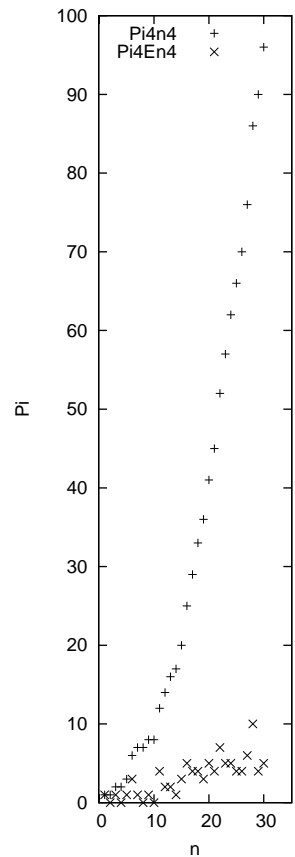}

\end{wrapfigure}%

Wegen der statistischen Unabhängigkeit aller Primzahlen voneinander
sollte man annehmen, daß auch die Vierlingsmenge unendlich ist. Dem
steht entgegen, daß es offenbar für kleinere Werte von n nur wenige
Quadratintervalle gibt, die einen Vierling enthalten. Da nun die Vierlinge
bezüglich ihrer Wahrscheinlichkeitsfunktion eine Menge der 4. Potenzordnung
darstellen, seien in \textbf{Abb. 8.1.} die Anzahl der Vierlinge $\leq n^{4}\ ,\ \pi_{4}\left(n^{4}\right),$
und die Vierlingsanzahl im Intervall $\eta_{4,n}:=\bigl(n^{4},\left(n+1\right)^{4}\bigr]\ ,\ \pi_{4}\left(\eta_{4,n}\right),$
über $n\leq30$ dargestellt. Offensichtlich ist das Intervall $\eta_{4,10}$
das größte vierlingsfreie Intervall, im Intervall $\eta_{4,14}$ tritt
letztmalig nur 1 Vierling auf. Diese Tatsache steht im Einklang mit
unserem statistischen Ansatz, nach dem zur Intervallbreite $b_{4,n}=4n^{3}+6n^{2}+4n+1$
\begin{equation}
\frac{\pi_{4}\left(\eta_{4,n}\right)}{b_{4,n}}=\left(\frac{1,06}{2\ln\left(n+1\right)}\right)^{4}\label{eq:8.1.}
\end{equation}
 als mittlere Primzahlvierlingsanzahl \\
$\leq\left(n+1\right)^{4}$ zu erwarten ist.

\textbf{Satz 8.1.} Da die Anzahl biquadrati\nobreakdash-scher Intervalle
mit mehr als 1 Vierling nach oben offen ist, muß es auch unendlich
viele Primzahlvierlinge geben.

Der Satz bleibt hier als Vermutung stehen. Ob es ein $n_{0}$ gibt,
sodaß $\forall n>n_{0}$ jedes Quadratintervall $\eta_{2,n}$ mindestens
einen Primzahlvierling enthält, ist damit noch nicht entschieden;
es ist aber nahe gelegt. Ein solches $n_{0}$ müßte aber deutlich
über $n=1000$ liegen, denn oberhalb $n=900$ gibt es noch primzahlvierlingsfreie
Quadratintervalle. Unterhalb $n=914$ gibt es insgesamt 96 Primzahlvierlinge.

\section{Zur Goldbach-Hypothese}

Christian Goldbach (1690-1764) vermutete in einem Brief an Leonhard
Euler, daß jede gerade Zahl $\geq6$ als Summe aus genau 2 ungeraden
Primzahlen dargestellt werden kann. Wir wollen sie in den etwas schärferen
Satz fassen:

\textbf{Satz 9.1.} a) Jede gerade Zahl $2m\geq8,\ m\epsilon\mathfrak{\mathbb{N}},$
kann als Summe aus genau 2 voneinander verschiedenen ungeraden Primzahlen
dargestellt werden. Zusätzlich \\
existiert $\forall$ Primzahlen p eine Darstellung $2m=2p.$ 

b) Die Darstellung $2m=p_{i}+p_{j},\ p_{i}\neq p_{j},$ ist im allgemeinen
vieldeutig; eindeutig ist sie nur für einige relativ kleine Werte
von $2m\leq2m_{0}=12$, (wenn $p_{i}=p_{j}$ mit gewertet wird). Die
Vielfachheit $v_{2m}$ der Goldbachpaar-Darstellungen von 2m erfüllt
\begin{equation}
v_{2m}=\left(n+\frac{1}{2}\right)\left(\frac{A^{'}}{\ln\left(n+1\right)}\right)^{2}\left(1+\delta_{2}\right),\ 2m\epsilon\eta_{n}.\label{eq:9.1.}
\end{equation}
 \textbf{Bew.:} Der Beweis soll auf probabilistischer Grundlage geführt
werden. Das ist zulässig, da die Gültigkeit der Verteilungsfunktion
$A^{'}/\ln x$ für jedes Intervall $\bigl(1,x\bigr]$ aus seiner Gültigkeit
$\forall\eta_{n}$ folgt, also bereits gezeigt ist, und alle Primzahlen
als statistisch unabhängig erwiesen sind. 

Die Summe aus 2 ungeraden Primzahlen 
\begin{equation}
p_{i}+p_{j}=2m\label{eq:9.2.}
\end{equation}
 ist stets geradzahlig. Aus der Forderung $p_{i}\neq p_{j}$ folgt
dann o.B.d.A. $p_{i}<p_{j}$ sowie die Existenz einer Zahl $\Delta$
, sodaß $p_{i}=m-\Delta$ und $p_{j}=m+\Delta$ gilt. Wir betrachten
nun ein Quadratintervall $\eta_{n}=\bigl(n^{2},\left(n+1\right)^{2}\bigr]$
mit $n\geq n_{0},$ sodaß auch noch jede gerade Zahl aus $\eta_{n-1}$
oberhalb m und das Intervall $\eta_{A,n}^{'}:=\bigl(0,4n\bigr]$ vollständig
unterhalb m liegt. Die Breite $b_{A}$ von $\eta_{A,n}^{'}$ ist gleich
der Summe der Breiten von $\eta_{n-1}$ und $\eta_{n}\ ,$ also $b_{A}=\left(2n+1\right)+\left(2n-1\right)=4n.$
Daher gilt $n_{0}=7.$ Für alle geraden Zahlen unterhalb dieser Größe
prüfen wir Satz 9.1. explizit numerisch:

4=2+2; 6=3+3; 8=3+5; 10=3+7=5+5; 12=5+7; 14=3+11=7+7; 16=3+13=5+11;
18=5+13=7+11; 20=3+17=7+13; 22=3+19=5+17=11+11; 24=5+19=7+17=11+13;
26=3+23=7+19=13+13; 28=5+23=11+17; 30=7+23=11+19; 32=3+29=13+19;\\
 34=3+31=5+29=11+23=17+17; 36=5+31=7+29=13+23=17+19.

Es sei hier bemerkt, daß in der Regel schon mindestens ein Goldbachpaar
auftritt mit $p_{j}\epsilon\eta_{n},$ in einigen Fällen wird aber
das erste Goldbachpaar erst mit $p_{j}\epsilon\eta_{n-1}$ gefunden.
Für sehr große n kann nicht ausgeschlossen werden, daß das erste $p_{j}$,
das ein Goldbachpaar bildet, noch kleiner ist. Das stört aber unsere
Betrachtung nicht, da zur Bestimmung der Vielfachheit v der möglichen
Goldbachpaar-Darstellungen alle $p_{j}$ mit $m\leq p_{j}<\left(n+1\right)^{2}$
zu berücksichtigen sind und hinreichend viele Quadratintervalle im
Intervall $\bigl(m,\left(n+1\right)^{2}\bigr]$ existieren.

Die Wahrscheinlichkeit $W_{g}$, daß die erste Primzahl in (\ref{eq:9.2.})
$p_{i}\epsilon\eta_{A,n}$ und die zweite $p_{j}\epsilon\eta_{n}$
erfüllt, ist durch das Produkt ihrer entsprechenden Wahrscheinlichkeiten
definiert. Wegen 
\begin{equation}
W_{i}=\frac{\pi\left(\eta_{A,n}\right)}{b_{n}}\approx2\frac{\pi\left(\eta_{n}\right)}{b_{n}}=2W_{j}\label{eq:9.3.}
\end{equation}
 gilt für das Paar $\left\{ p_{i},p_{j}\right\} $ 
\begin{equation}
W_{g}=2W_{j}^{2}.\label{eq:9.4.}
\end{equation}
 Für die einzelne Zahl 2m haben wir wieder den streuenden Anteil zu
verdoppeln, sodaß wir für die streuende Anzahl von Goldbachpaar-Darstellungen
zahlenmäßig das Doppelte der Primzahlzwillingsanzahl in $\eta_{n}$
gemäß (\ref{eq:7.9.}) erhalten: 
\begin{equation}
v_{2m}=\pi_{g}\left(2m\right)=2\pi_{2,s}\left(\eta_{n}\right)\ ,\label{eq:9.5.}
\end{equation}
 wie in Satz 9.1.b behauptet wurde, wobei $v_{2m}=\pi_{g}\left(2m\right)$
gesetzt ist;\textbf{ q.e.d.}

Beide Ausdrücke wachsen aber mit m über alle Schranken. Wegen des
Faktors 2 erhebt sich jedoch die Funktion $\pi_{g}\left(2m\right)$
inclusive ihrem Streubereich schneller über die Anzahl ``1'' . Die
\textbf{Abb. 9.1.} zeigt den Anfangsbereich von $\pi_{g}\left(2m\right)$
im numerischen Experiment zusammen mit den theoretisch ermittelten
Werten für $2m\leq330.$

\textbf{Satz 9.2.} Unter Goldbachpaaren versteht man im allgemeinen
nur Paare ungerader Primzahlen, deren Summe eine gerade Zahl 2m ergibt.
Läßt man auch die kleinste Primzahl $p_{1}=2$ als Paarpartner zu,
so können auch alle ungeraden Zahlen, die um 2 größer sind als eine
Primzahl, als Summe aus genau 2 Primzahlen dargestellt werden. Das
sind zwar unendlich viele ungerade Zahlen, genau so viele wie es Primzahlen
gibt, aber bei weitem nicht alle. Es gibt sogar unendlich viele Primzahlen,
die als Summe aus genau einer ungeraden Primzahl und $p_{1}=2$ dargestellt
werden können. Diese ungeraden Zahlen sind der kleinere Partner eines
jeden Primzahlzwillings, von denen wir zeigen konnten, daß ihre Anzahl
unendlich ist. 

\textbf{Satz 9.3.} Die Zahl 5 ist die einzige ungerade unberührbare
Zahl. \\
Man bezeichnet eine natürliche Zahl z als unberührbar, wenn es keine
natürliche Zahl x gibt, deren echte Teilersumme $\sigma^{*}\left(x\right)=z$
ist.

\textbf{Bew.:} Behauptung und Beweis sind angelehnt an Aufgabe 48
b in \cite{key-1}, p. 327 und 336. Die Behauptung folgt aus dem Goldbachtheorem
in der Fassung von Satz 9.1.a. Jede ungerade natürliche Zahl $z>8$
hat eine Darstellung $z=1+2n$ mit $2n\geq8$ . Jede gerade Zahl $2n\geq8$
hat mindestens eine Darstellung $2n=p+q$ mit den Primzahlen $p\neq q$
. Die Zahl $x=p\cdot q$ hat die Summe echter Teiler $\sigma^{*}\left(x\right)=1+p+q=1+2n=z.$
Die Zahlen $z=3$ und $z=7$ sind berührbar ($x=4$ bzw. $x=8$);
q.e.d.

\subsection*{Danksagung}

Für die Wartung und Systembetreuung meines PC danke ich Christian
Schmidt-Gütter und für Unterstützung bei der Arbeit mit \LaTeX{} und
gnuplot Susanne Gütter.
\end{document}